\documentclass[aap,seceqn,dvips]{arximspdf}

\doi{10.1214/105051606000000844}
\volume{17}
\issue{2}
\pubyear{2007}
\firstpage{572}
\lastpage{608}

\makeatletter
\newtheorem{theorem}{Theorem}[section]
\newproclaim{definition}[theorem]{Definition}
\newtheorem{lemma}[theorem]{Lemma}
\newtheorem{assume}[theorem]{Assumption}
\newtheorem{proposition}[theorem]{Proposition}
\newtheorem{corollary}[theorem]{Corollary}
\newproclaim{remark}[theorem]{Remark}
\newproclaim{example}[theorem]{Example}
\let\widebar\overline
\makeatother

\begin{document}
\begin{frontmatter}

\title{One-dimensional linear recursions with Markov-dependent coefficients}
\runtitle{Linear recursions with Markovian coefficients}

\begin{aug}
\author[A]{\fnms{Alexander} \snm{Roitershtein}\corref{}\ead[label=e1]{roiterst@math.ubc.ca}
\ead[label=u1,url]{www.math.ubc.ca/\texttildelow roiterst/}}
\runauthor{A. Roitershtein} \affiliation{University of British
Columbia}
\address[A]{Department of Mathematics\\
University of British Columbia\\
121-1984 Mathematics Road\\
Vancouver, British Columbia\\
Canada V6T 1Z2\\
\printead{e1}\\
\printead{u1}} %adresu isvedimo komanda gale!
\end{aug}

% HISTORY:
\received{\smonth{9} \syear{2004}} \revised{\smonth{10} \syear{2006}}

% ABSTRACT
%
\begin{abstract}
For a class of stationary Markov-dependent sequences $(A_n,B_n)
\in\mathbb{R}^2,$ we consider the random linear recursion $S_n=A_n+B_n
S_{n-1},$ $n \in\mathbb{Z},$ and show that the distribution tail of its
stationary solution has a power law decay.
\end{abstract}

% KEYWORDS
%

\end{frontmatter}

%s1 ###
\section{Introduction and statement of results}\label{s1}
Consider the stochastic difference equation
%
%e1.1 ###
\begin{equation}\label{main-def}
S_n=A_n+B_n S_{n-1},\qquad n \in{\mathbb N}, S_n \in{\mathbb R},
\end{equation}
with real-valued random coefficients $A_n$ and $B_n.$

If the sequence of random pairs $(A_n,B_n)_{n \in{\mathbb Z}}$ is
stationary and ergodic, $E (\log|B_0| )<0,$ and $E (\log|A_0|^+ )<
\infty,$ where $x^+=\max(0,x),$ then for any initial random value
$S_0,$ the limit law of $S_n$ is the same as that of the random
variable $R= A_0+\sum_{n=1}^{\infty} A_{-n} \prod_{i=0}^{n-1} B_{-i},$
and it is the unique initial distribution under which $(S_n)_{n \geq0}$
is stationary (cf. \cite{brandt}). Letting $\xi_n=A_{-n}$ and
$\rho_n=B_{-n}$ for $n \in{\mathbb Z},$ we get
%
%e1.2 ###
\begin{equation}\label{ar3} R=
\xi_0+\sum_{n=1}^{\infty} \xi_n \prod_{i=0}^{n-1} \rho_i.
\end{equation}
The stochastic difference equation \textup{(\ref{main-def})} has been
studied by many authors and has a remarkable variety of applications
(see, e.g., \cite{embre-goldie,samorachev,vervaat} for an extensive
account). The distribution tail of the random variable $R$ is the
topic of, for example,
\cite{goldie,goldie-grubel,grey,trakai75,kesten-randeq}, all assuming
that $(\rho_n,\xi_n)_{n \in{\mathbb Z}}$ is an i.i.d. sequence, and of
\cite{saporta}, where it is assumed that $(\rho_n)_{n \in{\mathbb Z}}$
is a finite Markov chain.

We study here the asymptotic behavior of the distribution tail of $R$
in the case that the sequence $(\zeta_n)_{n \in{\mathbb Z}}=
(\xi_n,\rho_n)_{n \in{\mathbb Z}}$ is an ``observable part'' of a
Markov-modulated process. By Markov-modulated process we mean the
following:

\begin{definition}
\label{case} Let $({\mathcal S},{\mathcal T})$ be a measurable space
and let $(x_n)_{n \in{\mathbb Z}}$ be a stationary Markov chain with
transition kernel $H(x,\cdot)$ defined on it.

A Markov-modulated process (MMP) associated with $(x_n)_{n \in{\mathbb
Z}}$ is a stationary Markov chain $(x_n,\zeta_n)_{n \in{\mathbb Z}}$
defined on a product space $({\mathcal S}\times\Upsilon, {\mathcal
T}\otimes\Xi ),$ whose transitions depend only on the position of
$(x_n).$ That is, for any $n \in{\mathbb Z}, A \in{\mathcal T}$,
\mbox{$B \in\Xi,$}
\[
P \bigl(x_n \in A, \zeta_n \in B | \sigma\bigl( (x_i,\zeta_i)\dvtx i<n
\bigr) \bigr)=\int_A H(x,d y) {\mathbb G}(x,y, B) |_{x=x_{n-1}},
\]
where ${\mathbb G}(x,y,\cdot)=P(\zeta_1 \in\cdot|x_0=x,x_1=y)$ is a
kernel on $({\mathcal S} \times{\mathcal S}\times\Xi).$
\end{definition}

For MMP $(x_n,\zeta_n)_{n \in{\mathbb Z}},$ where $\zeta_n=(\xi
_n,\rho_n),$ satisfying Assumption~\ref{measure1} below we show that
for some $\kappa>0,$ the limits $\lim_{t \to\infty} t^\kappa P(R>t)$
and $\lim_{t \to\infty} t^\kappa P(R<-t)$ exist and are not both zero.
Under our assumption, the parameter $\kappa$ is determined by
%
%e1.3 ###
\begin{equation}
\label{kappa239-43} \Lambda(\kappa)=0 ,\qquad\mbox{where }
\Lambda(\beta)=\lim_{n \to\infty} \frac{ 1}{ n } \log E
\Biggl(\,\prod_{i=0}^{n-1} |\rho_i|^\beta\Biggr).
\end{equation}
This extends both the one-dimensional version of a result of Kesten
which is valid for i.i.d. variables $(\xi_n,\rho_n)_{n \in{\mathbb Z}}$
(cf. \cite{kesten-randeq}, Theorem~5, see also an alternative approach
developed by Goldie in \cite{goldie}) as well as the recent result of
de Saporta \cite{saporta} where it is assumed that $(\rho_n)_{n
\in{\mathbb Z}}$ is a finite irreducible and aperiodic Markov chain
independent of the process $(\xi_n)_{n \in{\mathbb Z}}$ which is an
i.i.d. sequence.

In the joint paper with Eddy Mayer--Wolf and Ofer Zeitouni \cite{mars}, in the
context of an application to random walks in random environments, we
treated the particular case where $P(\xi_0=1,\rho_0>0)=1$ and
$(\rho_n)_{n \in{\mathbb Z}}$ is a point-wise transformation of a
stationary Markov chain $(x_n)_{n \in{\mathbb Z}}$ which is either
finite-state and irreducible (possibly periodic) or such that some
power of its transition kernel is dominated from above and below by a
probability measure (and thus is aperiodic).

The general case is more involved and requires additional arguments to
deal with it. We consider here the following Markov-modulated model
where the coefficients $(\xi_n,\rho_n)_{n \in{\mathbb Z}}$ of the
linear recursion \textup{(\ref{main-def})} are not necessarily
positive, the underlying Markov chain $(x_n)_{n \in{\mathbb Z}}$ is
defined in a general state space and may be periodic, and the sequences
$(\xi_n)_{n \in{\mathbb Z}}$ and $(\rho_n)_{n \in {\mathbb Z}}$ are not
assumed to be independent.

Let ${\mathcal B}$ denote the Borel $\sigma$-algebra of ${\mathbb R}.$

\begin{assume}
\label{measure1} There is a stationary Markov chain $(x_n)_{n \in
{\mathbb Z}}$ on a measurable space $ ({\mathcal S}, {\mathcal T})$
with transition kernel $H(x,\cdot)$ such that $(x_n,\zeta_n)_{n
\in{\mathbb Z}},$ where
\begin{eqnarray*}
\zeta_n=(\xi_n,\rho_n),\qquad n \in{\mathbb Z},
\end{eqnarray*}
is a MMP on $({\mathcal S}\times {\mathbb R}^2,{\mathcal
T}\times{\mathcal B}^{\otimes2})$ and:

\begin{longlist}[(A3)]
\item[(A1)] The $\sigma$-field ${\mathcal T}$ is countably generated.

\item[(A2)] The kernel $H(x,\cdot)$ is irreducible, that is, there
exists a $\sigma$-finite measure $\varphi$ on $({\mathcal S},{\mathcal
T})$ such that for all $x \in{\mathcal S},$ $\sum_{n=1}^\infty
H^n(x,A)>0$ whenever $\varphi(A)>0.$

\item[(A3)] There exist a probability measure $\mu$ on $({\mathcal
S},{\mathcal T}),$ a number $m_1 \in{\mathbb N},$ and a measurable
density kernel $h(x,y)\dvtx  {\mathcal S}^2 \to [0,\infty)$ such that
\[
H^{m_1}(x,A)=\int_A h(x,y) \mu(d y),
\]
and the family of functions $\{h(x,\cdot)\dvtx {\mathcal
S}\to[0,\infty)\} _{x \in {\mathcal S}}$ is uniformly integrable with
respect to the measure $\mu.$

\item[(A4)] $P(|\xi_0| < c_\xi)=1$ for some $c_\xi>0.$

\item[(A5)] $P(c_\rho^{-1}< |\rho_0| <c_\rho)=1$ for some $c_\rho>1.$

\item[(A6)] Let $\Lambda(\beta)=\limsup_{n \to\infty} \frac{ 1}{ n }
\log E (\prod_{i=0}^{n-1} |\rho_i|^{\beta_1} ).$ Then there exist
constants $\beta_1>0$ and $\beta_2>0$ such that $\Lambda(\beta _1)\geq
0$ and $\Lambda(\beta_2)<0.$

\item[(A7)] There do not exist a constant $\alpha>0$ and a measurable
function $\beta\dvtx {\mathcal S}\times\{ -1,1\} \to [0,\alpha)$ such
that
\[
P \bigl(\log|\rho_1| \in\beta(x_0,\eta)-
\beta\bigl(x_1,\eta\cdot\operatorname{sign} (\rho_1)
\bigr)+\alpha\cdot{\mathbb Z} \bigr)=1,
\]
for $\eta\in\{-1,1\}.$
\end{longlist}
\end{assume}

\begin{remark}
\label{remr} The assumption that the sequence $(x_n,\zeta_n)_{n
\in{\mathbb Z}}$ is stationary is explicit in Definition~\ref{case} of
Markov-modulated processes. It turns out (see Lemma~\ref{urg3} below)
that under assumptions (A1)--(A3), the Markov chain $(x_n)_{n
\in{\mathbb Z}}$ has a unique stationary distribution. This
distribution induces a (unique) stationary probability measure for the
sequence (Markov chain) $(x_n,\zeta_n)_{n \in{\mathbb Z}},$ which we
denote by $P.$ The expectation according to the stationary measure $P$
is denoted by~$E.$
\end{remark}

Note that condition (A6) implies by Jensen's inequality that
$E(\log|\rho_0|)<0.$ Thus, by a theorem of Brandt \cite{brandt}, the
series in \textup{(\ref{ar3})} converges absolutely, $P\mbox{-a.s.}$ It
will be shown later (see Lemma~\ref{urg5} below) that the both
$\limsup$'s in (A6) is in fact a limit, and thus this condition
guarantees, by convexity, the existence of a unique $\kappa$ in
\textup{(\ref{kappa239-43})}.

Assumption (A7) ensures that $\log|\rho_n|$ is nonarithmetic (in the
sense of the following definition) relative to both the underlying
process $(x_n)_{n \in{\mathbb Z}}$ as well as to the auxiliary chain
$(\hat x_n)_{n \in{\mathbb Z}}$ introduced in Section~\ref{aux}.

\begin{definition}[(\cite{alsmeyer,shur})]\label{arc}
Let $(x_n,q_n)_{n \in{\mathbb Z}}$ be a MMP. The
process $(q_n)_{n \in{\mathbb Z}}$ is said to be nonarithmetic relative
to the Markov chain $(x_n)_{n \in{\mathbb Z}}$ if there do not exist a
constant $\alpha
>0$ and a measurable function $\beta\dvtx {\mathcal S}\to[0,\alpha)$ such that
\[
P \bigl(q_0 \in\beta(x_{-1})-\beta(x_0)+\alpha\cdot{\mathbb Z} \bigr)=1.
\]
\end{definition}

We will next state our results for the coefficients $(\xi_n,\rho_n)_{n
\in{\mathbb Z}}$ satisfying Assumption~\ref {measure1}. We will denote
%
%e1.4 ###
\begin{equation}\label{cond43}
P_x^{-}(\cdot)=P(\cdot|x_{-1}=x) \quad\mbox{and}\quad
E_x^{-}(\cdot)=E(\cdot|x_{-1}=x),
\end{equation}
and keep the notation $P_x(\cdot)$ and $E_x(\cdot)$ for
$P(\cdot|x_0=x)$ and $E(\cdot|x_0=x)$, respectively.

The case of positive coefficients $(\xi_n,\rho_n)_{n \in{\mathbb Z}}$
is qualitatively different from and technically simpler than the
general one [e.g., it turns out that in this case $\lim_{t \to\infty}
t^\kappa P(R>t)$ is \textit{always} positive], and it will be
convenient to treat it separately.

\begin{theorem} \label{main-rec}
Let Assumption~\textup{\ref{measure1}} hold and denote by $\pi$ the
stationary distribution of the Markov chain $(x_n)_{n \in{\mathbb Z}}.$

If $P(\xi_0>0,\rho_0>0)=1$ then for $\pi$-almost every $x \in {\mathcal
S},$ the following limit exists and is strictly positive:
\begin{eqnarray*}
K(x)=\lim_{t \to\infty} t^\kappa P_x^{-}( R
>t),
\end{eqnarray*}
where the parameter $\kappa$ is given by \textup{(\ref{kappa239-43})}
and the random variable $R$ is defined in~\textup{(\ref{ar3})}.
\end{theorem}

An application of Theorem~\ref{main-rec} and estimates
\textup{(\ref{bound})}, \textup{(\ref{43-bound-43})} to random walks in
random environments can be found in \cite{mars}. The main step of the
proof follows Goldie's argument (cf. \cite{goldie}, Theorem~2.3)
closely and relies on the application of a version (due to Alsmeyer,
cf. \cite{alsmeyer}) of the Markov renewal theorem due to Kesten (cf.
\cite{kesten-rem}, see also \cite{mac-dan,shur} and references to
related articles in \cite{kesten-rem}).

For coefficients $(\xi_n,\rho_n)_{n \in{\mathbb Z}}$ with arbitrary
signs we have:

\begin{theorem}
\label{main-rec-3} Let Assumption~\textup{\ref{measure1}} hold and
denote by $\pi$ the stationary distribution of the Markov chain
$(x_n)_{n \in{\mathbb Z}}.$

Then, with $\kappa$ given by \textup{(\ref{kappa239-43})} and $R$
defined in \textup{(\ref{ar3})},
\begin{longlist}[(a)]
\item[(a)] For $\pi$-almost every $x \in{\mathcal S},$ the following
limits exist:
%
%e1.5 ###
\begin{equation}
\label{gvul1} K_1(x)=\lim_{t \to\infty} t^\kappa P_x( R
>t)\quad\mbox{and}\quad K_{-1}(x)=\lim_{t \to\infty} t^\kappa
P_x(R<-t).
\end{equation}

\item[(b)] $\pi(K_1(x)+K_{-1}(x)>0 )\in\{0,1\}.$

\item[(c)] If Condition \textup{G} (see
Definition~\textup{\ref{congee}} below) is satisfied then it holds that
$\pi(K_1(x)=K_{-1}(x) )=1.$ Moreover, if Condition \textup{G} is not
satisfied then either $\pi(K_1(x)>0 \mbox{ and } K_{-1}(x)>0 )\in\{0,1\}$
or there exists a (possibly trivial) partition of ${\mathcal S}$ into
two disjoint measurable sets $A$ and $B$ such that $\pi$-a.s.,
$K_1(x)>0$ and $K_{-1}(x)=0$ for $x \in A$ whereas $K_1(x)=0$ and
$K_{-1}(x)>0$ for $x \in B.$
\end{longlist}
\end{theorem}

\begin{definition}
\label{congee} We say that Condition \textup{G} holds if there does not
exist a (possibly trivial) partition of ${\mathcal S}$ into two
disjoint measurable sets $A_1$ and $A_{-1}$ such that for $i
\in\{-1,1\},$
\[
P(x_0 \in A_i, x_1 \in A_{-i},\rho_1>0)= P(x_0 \in A_i, x_1 \in
A_i,\rho_1<0)=0.
\]
\end{definition}

Condition G is a generalization of the condition of $l$-irreducibility
introduced in \cite{saporta}. Note that this condition is not satisfied
if $P(\rho_0>0)$ (take $A_1={\mathcal S}$ and $A_{-1}=\varnothing$).
Proposition~\ref{still} shows that Condition G is equivalent to the
assertion that the Markov chain $\hat x_n=(x_n,\gamma_n),$ where
$\gamma_n=\operatorname{sign}(\rho_0 \cdots\rho_{n-1}),$ is irreducible
under Assumption~\ref{measure1}.

The proof of Theorem~\ref{main-rec-3} is basically by applying a
Markovian adaptation of the implicit renewal theory of Goldie
\cite{goldie} (see Section~\ref{proofs}) to the Markov chain $\hat x_n$
and the random walk $V_n=\sum_{i=0}^{n-1} \log|\rho_i|.$ The Markov
chain $\hat x_n$ carries the necessary information about the sign of
the products of $\rho_i$ and at the same time, as we shall see in
Section~\ref{aux}, inherits all essential properties of the Markov
chain~$x_n.$

In order to show that $K_1(x)+K_{-1}(x)>0$ in Theorem \ref{main-rec-3},
we need an extra nondegeneracy assumption which guarantees that the
random variable $R$ is not a deterministic function of the initial
state $x_{-1}.$ Again following \cite{goldie} and using the renewal
theory developed in \cite{alsmeyer}, we complement Theorem
\ref{main-rec-3} by the following necessary and sufficient condition
for $R$ to be nondegenerate under $P_x^{-}$ and for the limit to be
positive. This condition is a natural generalization of the criterion
that appears in the case where the random variables $(\xi_n,\rho_n)$
are i.i.d. (cf. \cite{kesten-randeq} and \cite{goldie}). Note that the
condition is trivially satisfied under the assumptions of
\cite{saporta} [because $(\xi_n)_{n \in {\mathbb Z}}$ is assumed to be
independent of $(x_n)_{n \in{\mathbb Z}}$].

\begin{theorem} \label{mishne}
Let Assumption~\textup{\ref{measure1}} hold and denote by $\pi$ the
stationary distribution of the Markov chain $(x_n)_{n \in{\mathbb Z}}.$
Then:
\begin{longlist}[(a)]
\item[(a)] $\pi(K_1(x)+K_{-1}(x)>0 )=0$ if and only if there exists a
measurable function $\Gamma\dvtx  {\mathcal S}\to{\mathbb R}$ such that
%
%e1.6 ###
\begin{equation}\label{aper}
P \bigl(\xi_0+\Gamma(x_0)\rho_0=\Gamma(x_{-1}) \bigr)=1.
\end{equation}

\item[(b)] There exists a constant $C_1>0$ such that for $\pi$-almost
every $x \in{\mathcal S},$
%
%e1.7 ###
\begin{equation}\label{bound} t^\kappa
P^{-}_x ( |R|>t )\leq C_1\qquad\forall t>0.
\end{equation}
In particular, $\lim_{t \to\infty} t^\kappa P(R>t)=E (K_1(x_{-1}) ),$
$\lim_{t \to\infty} t^\kappa P(R<-t)=\break E (K_{-1}(x_{-1}) ),$ and the
limits are finite.

\item[(c)] If \textup{(\ref{aper})} does not hold for any measurable
function $\Gamma\dvtx  {\mathcal S}\to{\mathbb R},$ then there exist
positive constants $C_2$ and $t_c$ such that for $\pi$-almost every $x
\in{\mathcal S},$
%
%e1.8 ###
\begin{equation}
\label{43-bound-43} t^\kappa P^{-}_x ( |R|>t )\geq C_2\qquad \forall
t>t_c.
\end{equation}
In particular, $\lim_{t \to\infty} t^\kappa P( R>t)$ and $\lim_{t
\to\infty} t^\kappa P(R<-t)$ are not both zero.
\end{longlist}
\end{theorem}

\begin{remark}
Throughout this paper we work with the probability measures
$P^{-}_x(\cdot)=P(\cdot|x_{-1}=x)$ defined in \textup{(\ref{cond43})}
rather than with $P_x(\cdot)=P(\cdot|x_0=x).$ Since
\[
P_x(R>t)=E \biggl( P^{-}_x \biggl(R>\frac{ t-a}{ b } \biggr) \Big|
\xi_0=a,\rho_0=b \biggr),
\]
the bounded convergence theorem and part (b) of
Theorem~\ref{main-rec-3} show that all our results hold also for the
usual conditional measure $P_x.$

However, treating the linear recursion \textup{(\ref{main-def})} in the
setup of Markov-modulated processes, it is not so natural to work with
the conditional probabilities $P_x.$ In order to elucidate this point,
let us consider the following two examples:
\begin{longlist}[(ii)]
\item[(i)] The random variable $R$ is conditionally independent of the
\textup{``}past,\textup{''} that is, of $\sigma ((\xi_n,\rho_n)_{n<0} ),$
given $x_{-1}$ but not given $x_0.$

\item[(ii)] Let $\tau>0$ be a finite random time such that $x_\tau$ is
distributed according to a probability measure $\psi,$ and define
$P_\psi(\cdot):= \int_{\mathcal S}P_x(\cdot)\psi(dx)$ and
$R_\tau=\xi_\tau+\sum_{n=\tau+1}^\infty\xi_n
\prod_{i=\tau}^{n-1}\rho_i.$ Then in general, since the distribution of
$x_{\tau-1}$ and hence that of $(\xi_\tau,\rho_\tau)$ are unknown,
\[
P(R_\tau\in\cdot) \neq P_\psi(R \in\cdot).
\]
On the other hand, $P(R_{\tau+1} \in\cdot)=P_\psi^{-}(R \in\cdot)$ for
$P_\psi^{-}(\cdot):= \int_{\mathcal S}P_x^{-}(\cdot)\psi(dx).$
\end{longlist}
\end{remark}

The rest of the paper is organized as follows. Section \ref{proper},
divided into three subsections, is mostly devoted to the properties of
the Markov chain $(x_n)_{n \in{\mathbb Z}}$ and of the random walk
$V_n=\sum_{i=0}^{n-1} \log|\rho_i|.$ Section~\ref{bagr} is devoted to
the basic properties of the underlying Markov chain $(x_n)_{n
\in{\mathbb Z}}.$ In Section~\ref{pft} we state a Perron--Frobenius
type theorem (Proposition~\ref{urg}) which plays an important role in
the subsequent proofs and in particular implies the existence and
uniqueness of $\kappa$ in \textup{(\ref{kappa239-43})} (see Lemma
\ref{urg5}). The proof of Proposition~\ref{urg} is deferred to the
\hyperref[proofp]{Appendix}. Section~\ref{similar} is devoted to the
Markov renewal theory which is then used in Section~\ref{arbitar},
where it is applied to the Markov chain $\hat x_n=(x_n,\gamma_n)$ and
the random walk $V_n.$ Section~\ref{proofs} contains a reduction of
Theorems \ref{main-rec} and~\ref{main-rec-3} to a renewal theorem which
is an adaptation of a particular case of Goldie's implicit renewal
theorem (cf. \cite{goldie}). Section~\ref{aux} is devoted to study of
the auxiliary Markov chain $\hat x_n.$ The main goal here is to show
that the renewal theorem obtained in Section~\ref{proofs} can be
applied to the couple $(\hat x_n,V_n).$ The proofs of the main results
(Theorems \ref{main-rec},~\ref{main-rec-3} and~\ref{mishne}) are then
completed in Section~\ref{arbitar}.

%s2 ###
\section{Background and preliminaries}
\label{proper} Similarly to the i.i.d. case (cf. \cite{goldie}
and \cite{kesten-randeq}), the asymptotic behavior of the tail of $R$
under Assumption~\ref{measure1} is determined by the properties of
$V_n= \sum_{i=0}^{n-1} \log|\rho_i|$ and in particular is closely
related to the renewal theory for this random walk. This section is
devoted to the properties of the Markov chain $(x_n)_{n \in{\mathbb
Z}}$ and of the associated random walk with Markov-dependent
increments. The aim here is to provide for future use some technical
tools, namely the regeneration times $N_i$ defined in
Section~\ref{bagr} by the Athreya--Ney--Nummelin procedure, a
Perron--Frobenius theorem for positive kernels stated in
Section~\ref{pft}, and the Markov renewal theory recalled in
Section~\ref{similar}.

%s2.1 ###
\subsection{Some properties of the underlying Markov chain $(x_n)_{n
\in{\mathbb Z}}$} \label{bagr} First, let us note that assumption (A3)
implies that the transition kernel $H$ is \textit{quasi-compact}.
Recall that a transition probability kernel $H(x,\cdot)$ on a
measurable space $({\mathcal S},{\mathcal T})$ is called quasi-compact
if there exist constants $\varepsilon\in(0,1),$ $\delta\in(0,1),$ $m_1
\in{\mathbb N},$ and a probability measure $\mu$ such that
$H^{m_1}(x,A)<\varepsilon$ whenever $\mu(A)<\delta,$ or alternatively,
$H^{m_1}(x,A)>1-\varepsilon$ whenever $\mu(A)>1-\delta.$ If a
quasi-compact kernel $H$ is the transition kernel of a Markov chain
$(x_n)_{n \in{\mathbb Z}},$ then the chain is also called
quasi-compact. The condition on transition kernels used in this
definition was introduced by Doeblin (see, e.g., \cite{yosida} for a
historical account).

In the following lemma we summarize some properties of quasi-compact
chains which will be useful in the sequel (see Theorem 3.7 in
\cite{revuz}, Chapter~6, Section~3 for the first three assertions,
Proposition 5.4.6 and Theorem 16.0.2 in \cite{meyn-tweedie} for the
fourth, and Propositions 3.5, 3.6 in \cite{revuz}, Chapter~3, Section~3
for the last one).

\begin{lemma} \label{urg3}
Let $(x_n)_{n \in{\mathbb Z}}$ be an irreducible quasi-compact Markov
chain defined on a measurable space $({\mathcal S},{\mathcal T}).$
Then, there exist a number $d \in {\mathbb N}$ \textup{[}the period of
$(x_n)_{n \in{\mathbb Z}},$\textup{]} a sequence of $d$ disjoint measurable sets
$({\mathcal S}_1,{\mathcal S}_2,\ldots,{\mathcal S}_d)$ (a $d$-cycle),
and probability measures $\pi$ and $\psi$ on $({\mathcal S},{\mathcal
T})$ such that:
\begin{longlist}[(iii)]
\item[(i)] The following holds for all $i=1,\ldots,d,$ and $x
\in{\mathcal S}_i$: $H(x,{\mathcal S}_j^c)=0$ for $j=i+1$
$(\operatorname{mod} d).$

\item[(ii)] $\pi$ is the unique stationary distribution of $(x_n),$
$\pi({\mathcal S}_i)>0$ for $i=1,\ldots,d,$ and $\pi({\mathcal
S}_0)=1,$ where ${\mathcal S}_0=\bigcup_{i=1}^d {\mathcal S}_i.$

\item [(iii)] $(x_n)_{n \in{\mathbb Z}}$ is Harris recurrent chain when
restricted to the states of the set~${\mathcal S}_0.$ That is, $P(x_n
\in A \mbox{ i.o. for} n \geq0|x_0=x)=1,$ for all $x \in{\mathcal
S}_0$ and measurable $A \subseteq{\mathcal S}_0$ with $\pi(A)>0.$

\item[(iv)] $\psi({\mathcal S}_1)=1,$ and there exist constants $r
\in(0,1)$ and $m \in{\mathbb N}$ such that
%
%e2.1 ###
\begin{equation}\label{lem-lem} H^m(x,A) > r \psi(A)\qquad
\forall  x \in{\mathcal S}_1, A \in{\mathcal T}.
\end{equation}

\item[(v)] The process $(x_n)_{n \in{\mathbb Z}}$ is ergodic under its
stationary distribution.
\end{longlist}
\end{lemma}

The minorization condition \textup{(\ref{lem-lem})} with \textit{ some}
recurrent set ${\mathcal S}_1$ is equivalent to the Harris recurrence
(see, e.g., \cite{nummelin}). The particular form of the set ${\mathcal
S}_1$ in (iv) as cyclic element is particularly advantageous and is due
to the Doeblin condition.

We will next define a sequence of regeneration times $\{N_i\}_{i \geq
0}$ for the Markov chain $(x_n)_{n \in{\mathbb Z}}$ restricted to
$({\mathcal S}_0,{\mathcal T}_0),$ where ${\mathcal T}_0=\{A \in
{\mathcal T}\dvtx  A \subseteq {\mathcal S}_0\}.$ Let the set
${\mathcal S}_1$ and the number $m$ be the same as in
\textup{(\ref{lem-lem})}, and let $N_0$ be the first hitting time of
the set ${\mathcal S}_1$:
%
%e2.2 ###
\begin{equation}\label{nnn}
N_0=\inf\{n \geq-1 \dvtx  x_n \in{\mathcal S}_1\}.
\end{equation}
Note that $N_0 \leq d-1$ and $N_0$ is a deterministic function of
$x_{-1}$ on the set ${\mathcal S}_0.$ The randomized stopping times
$N_i,$ $i \geq1,$ can be defined in an enlarged (if needed) probability
space by the following procedure (see \cite{asmus,athreya-ney,nummelin}).
Given a state $x_{{N_0}} \in{\mathcal S}_1,$ generate
$x_{{N_0+m}}$ as follows: with probability $r$ distribute $x_{{N_0+m}}$
over ${\mathcal S}_0$ according to $\psi$ and with probability $1-r$
according to $1/(1-r)\cdot\Theta(x_0,\cdot),$ where the substochastic
kernel $\Theta(x,\cdot)$ is defined by
%
%e2.3 ###
\begin{equation}
\label{theta} H^m(x,A)=\Theta(x,A)+r \mathbf{1}_{{\mathcal
S}_1}(x)\psi(A),\qquad x \in {\mathcal S}_0, A \in{\mathcal T}_0.
\end{equation}
Then, (unless $m=1$) sample the segment $
(x_{{N_0+1}},x_{{N_0+2}},\ldots, x_{{ N_0+m-1}} )$ according to the
$(x_n)_{n \in{\mathbb Z}}$ chain's conditional distribution, given
$x_{{N_0}}$ and $x_{{ N_0+m}}.$ Generate $x_{{N_0+2m}}$ and
$x_{{N_0+m+1}},$ $x_{{N_0+m+2}}, \ldots, x_{{N_0+2m-1}}$ in a similar
way, and so on. Let $\{n_j\}_{j \geq1}$ be the successful times when
the move of the chain $(x_{{N_0+m n}})_{n \geq0}$ is according to
$\psi,$ and set $N_j=N_0+m n_j,$ $j \geq1.$ Note that $N_j$ is not the
$j$th visit to~${\mathcal S}_1.$

By construction, the blocks $ (x_{{ N_{_i}+1}},x_{{N_{_i}+2}}, \ldots,
x_{{N_{{i+1}}}} )$ are one-dependent and for $i \geq1$ they are
identically distributed ($x_{{N_{_i}}},$ $i \geq1,$ are independent and
distributed according to the measure $\psi$). It follows from the
construction that the random times $N_{i+1}-N_i$ are i.i.d. for $i
\geq0,$ and that there exist constants $\vartheta\in {\mathbb N},$
$\delta>0,$ such that
%
%e2.4 ###
\begin{equation}\label{as} P^{-}_x(N_1 \leq
\vartheta)> \delta\qquad \forall x \in{\mathcal S}_0.
\end{equation}
We summarize the properties of the random times $N_i$ in the following
lemma.

\begin{lemma} \label{ndef433}
Let $(x_n)_{n \in{\mathbb Z}}$ be an irreducible quasi-compact Markov
chain with state space ${\mathcal S},$ and let the set ${\mathcal S}_0$
be as in Lemma~\textup{\ref{urg3}}.

Then there exists a strictly increasing sequence $(N_i)_{i\ge0}$ of
random times such that:
\begin{longlist}[(iii)]
\item[(i)] $(N_{i+1}-N_i)_{i\ge0}$ are i.i.d.

\item[(ii)] The blocks $(x_{{N_i+1}}^{(0)},\ldots,x_{{N_{i+1}}}^{(0)})$
are one-dependent for $i \geq0$ and identically distributed for
$i\ge1$, where $(x_n^{(0)})_{n \in{\mathbb Z}}$ is the Markov chain
induced by $(x_n)_{n \in{\mathbb Z}}$ on $({\mathcal S}_0,{\mathcal
T}_0).$

\item[(iii)] $N_0 \leq d-1, \forall x \in{\mathcal S}_0,$ where $d$ is
the period of $(x_n)_{n \in{\mathbb Z}}.$
\item[(iv)] There exist constants $\vartheta\in{\mathbb N}$ and
$\delta>0$ such that \textup{(\ref{as})} is satisfied.
\end{longlist}
\end{lemma}

Throughout the rest of the paper we shall be concerned with the
measurable space $({\mathcal S}_0,{\mathcal T}_0),$ where ${\mathcal
T}_0=\{A \in{\mathcal T}\dvtx   A \subseteq{\mathcal S}_0\},$ rather
than with $({\mathcal S},{\mathcal T}).$ Without loss of generality we
may and shall assume that
%
%e2.5 ###
\begin{equation}\label{calso}
P_x^{-} \bigl(|\xi_0| < c_\xi \mbox{ and } |\rho_0|
\in(c_\rho,^{-1},c_\rho) \bigr)=1\qquad \forall x \in{\mathcal S}_0.
\end{equation}
Otherwise we can restrict our attention to the Markov chain induced by
$(x_n)_{n \in{\mathbb Z}}$ on the set of full measure $\pi$ where the
equality in \textup{(\ref{calso})} does hold. Clearly,
Assumption~\ref{measure1} and Lemma~\ref{urg3} remain true for this
Markov chain.

%s2.2 ###
\subsection{A Perron--Frobenius theorem for positive bounded kernels}
\label{pft} The aim of this subsection is to state a Perron--Frobenius
theorem for positive kernels (Proposition~\ref{urg} below). Proposition
\ref{urg} is an essential part of the subsequent proofs where it is
applied to kernels of the form $K(x,A)=E_x^{-} (\prod_{i=0}^n
|\rho_i|^\beta; x_n \in A )$ and $\widehat \Theta(x,A)=E_x^{-}
(\mathbf{1}_{\{n <N_1\}}\prod_{i=0}^n |\rho _i|^\beta; x_n \in A ),$
where the random time $N_1$ is defined in Section~\ref{bagr}. The proof
of Proposition~\ref{urg} is deferred to Appendix~\ref{proofp}.

One immediate consequence of this proposition is the following lemma
which proves the existence and uniqueness of the parameter $\kappa$ in
\textup{(\ref{kappa239-43})}.

\begin{lemma}
\label{urg5} Let Assumption~\textup{\ref{measure1}} hold and let the
set ${\mathcal S}_0$ be as defined in Lemma~\textup{\ref{urg3}}. Then,
\begin{longlist}[(a)]
\item[(a)] For any $\beta>0$ and every $x \in{\mathcal S}_0,$ the
following limit exists and does not depend on $x$:
%
%e2.6 ###
\begin{equation}
\label{kappa239} \Lambda(\beta)=\lim_{n \to\infty} \frac{ 1}{ n } \log
E_x^{-} \Biggl(\,\prod_{i=0}^{n-1} |\rho_i|^\beta\Biggr).
\end{equation}
Moreover, for some constants $c_\beta\geq1$ that depend on $\beta$
only,
%
%e2.7 ###
\begin{equation}\label{h-betaa} c_\beta^{-1} e^{n \Lambda(\beta)}
\leq E_x^{-} \Biggl(\,\prod_{i=0}^{n-1} |\rho_i|^\beta\Biggr) \leq
c_\beta e^{n \Lambda(\beta)}\qquad \forall x \in{\mathcal S}_0, n
\in{\mathbb N}.
\end{equation}

\item[(b)] There exists a unique $\kappa>0$ such that
$\Lambda(\kappa)=0,$ $\Lambda(\beta)(\beta-\kappa) \geq0$ for all
$\beta>0.$
\end{longlist}
\end{lemma}

We next proceed with Proposition~\ref{urg}, from which the lemma is
derived at the end of this subsection.

A function $K\dvtx  {\mathcal S}_0 \times{\mathcal T}_0 \to(0,\infty)$
is a \textit{ positive bounded kernel}, or simply \textit{ kernel}, if
the following three conditions hold: (i) $K(\cdot,A)$ is a measurable
function on ${\mathcal S}_0$ for all $A \in{\mathcal T}_0,$ (ii)
$K(x,\cdot)$ is a finite positive measure on ${\mathcal T}_0$ for all
$x \in{\mathcal S}_0,$ (iii) $\sup_{x \in{\mathcal S}_0} K(x,{\mathcal
S}_0) <\infty.$ Let $B_b$ be the Banach space of bounded measurable
real-valued functions on the measurable space $({\mathcal
S}_0,{\mathcal T}_0)$ with the norm $\|f\|=\sup_{x \in{\mathcal
S}_0}|f(x)|.$ Any positive bounded kernel $K(x,A)$ defines a bounded
linear operator on $B_b$ by setting $K f(x)= \int_{{\mathcal S}_0}
K(x,d y)f(y).$ We denote by $r_{_K}$ the spectral radius of the
operator corresponding to the kernel $K,$ that is
\[
r_{_K}=\lim_{n \to\infty} \sqrt[n]{\rule{0pt}{7pt}\smash{\|K^n
\mathbf{1}\|}}=\lim_{n\to\infty} \sqrt[n]{\|K^n\|},
\]
where $\mathbf{1}(x)\equiv1.$

The following proposition generalizes Lemma 2.6 in \cite{mars} allowing
us to deal with a more general class of underlying Markov chains
$(x_n)_{n \in{\mathbb Z}}.$

\begin{proposition}
\label{urg} Let $K(x,\cdot)$ be a positive bounded kernel on
$({\mathcal S}_0,{\mathcal T}_0)$ and $s(x,y)\dvtx  {\mathcal S}^2_0
\to {\mathbb R}$ be a measurable function such that $s(x,y)
\in(c_1^{-1},c_1)$ for some $c_1>1$ and all $(x,y) \in{\mathcal
S}^2_0.$ Assume that there exists a set ${\mathcal S}_1 \in{\mathcal
T}_0$ such that:
\begin{longlist}[(iii)]
\item[(i)] For some constants $d \in{\mathbb N},p>0,$
\[
\sum_{i=1}^d K^i(x,{\mathcal S}_1) \geq p\qquad \forall  x \in{\mathcal
S}_0.
\]

\item[(ii)] For some constant $m \in{\mathbb N}$ and probability
measure $\psi$ concentrated on ${\mathcal S}_1,$
\[
K^m(x,{\mathcal S}_1^c)=0\qquad \forall x \in{\mathcal S}_1,
\]
where ${\mathcal S}_1^c$ is the complement set of ${\mathcal S}_1,$ and
%
%e2.8 ###
\begin{equation}
\label{comp1} K^m(x,A) \geq\int_A s(x,y) \psi(d y)\qquad \forall  x
\in{\mathcal S}_1, A \in{\mathcal T}_0.
\end{equation}
Further, assume that:

\item[(iii)] There are a probability measure $\mu$ on $({\mathcal
S}_0,{\mathcal T}_0)$ and a constant $m_1 \in{\mathbb N}$ such that for
all $\varepsilon>0$ there exists $\delta=\delta(\varepsilon)>0$ such
that
%
%e2.9 ###
\begin{equation}
\label{doeblin} \mu(A)<\delta\qquad\mbox{implies } \sup_{x \in
{\mathcal S}_0} K^{m_1}(x,A) <\varepsilon.
\end{equation}
\end{longlist}

\textup{[}This condition entails $K(x,\cdot)\ll\mu$ for all $x
\in{\mathcal S}_0.$\textup{]}

Let ${\mathcal T}_1=\{A \in {\mathcal T}_0\dvtx  A \subseteq{\mathcal
S}_1\}$ and let a kernel $\widehat\Theta(x,A)$ on $({\mathcal
S}_1,{\mathcal T}_1)$ be such that
\[
K^m(x,A)=\widehat \Theta(x,A)+r \int_A s(x,y) \psi(d y)\qquad \forall x
\in{\mathcal S}_1,A \in {\mathcal T}_1,
\]
for some $r \in(0,1).$

Then:
\begin{longlist}[(a)] \item[(a)] There exists a function $f\in
B_b$ such that $ \inf_x f(x)>0$ and $K f=r_{_K} f.$

\item[(b)] There exists a function $g \in B_b$ such that $ \inf_x
g(x)>0$ and $\widehat\Theta g=r_{\widehat\Theta} g.$

\item[(c)] $r_{\widehat\Theta}\in(0,r_{_K}^m).$
\end{longlist}
\end{proposition}

The proof of the proposition is included in Appendix~\ref{proofp}.

\begin{pf*}{Proof of Lemma~\ref{urg5}}\label{harris1}
Let $Q(x,y,B)=P(\rho_{-n} \in B|x_{n-1}=x,x_n=y),$ and for any
$\beta\geq0$ define the kernel $H_\beta(x,\cdot)$ on $({\mathcal
S}_0,{\mathcal T}_0)$ by
%
%e2.10 ###
\begin{equation}\label{kbeta} H_\beta(x,dy)=H(x,dy)
\int_{\mathbb R}Q(x,y,dz) |z|^\beta.
\end{equation}
Then for any $\beta\geq0,$
%
%e2.11 ###
\begin{equation}\label{hbeta} E_x^{-}
\Biggl(\,\prod_{i=0}^{n-1} |\rho_i|^\beta\Biggr)= H_\beta^n \mathbf{1}
(x)\qquad \forall x \in{\mathcal S}_0.
\end{equation}
The kernels $H_\beta, \beta\geq0,$ satisfy the conditions of
Proposition~\ref{urg}. It follows from \textup{(\ref{bbb})} with
$K=H_\beta$ that for some constant $c_\beta\geq1$ which depends on
$\beta$ only,
%
%e2.12 ###
\begin{equation}\label{h-beta} c_\beta^{-1}
r_\beta^n \leq E_x^{-} \Biggl(\,\prod_{i=0}^{n-1} |\rho_i|^\beta\Biggr)
\leq c_\beta r_\beta^n\qquad \forall x \in{\mathcal S}_0, n \in{\mathbb
N},
\end{equation}
where $r_\beta=r_{H_\beta}.$ This yields assertion (a) of the lemma.
The claim of its part (b) follows then from the convexity of the
function $\Lambda(\beta)$ which takes by assumption~(A6) both positive
and negative values.
\end{pf*}

%s2.3 ###
\subsection{Markov renewal theory}
\label{similar} The proofs of our results rely on the use of the
following version of the Markov renewal theorem which is due to
Alsmeyer \cite{alsmeyer}. Recall Definition~\ref{case} of
Markov-modulated processes and Definition~\ref{arc} of nonarithmetic
processes. Let ${\mathcal B}$ denote the Borel $\sigma$-algebra of
${\mathbb R}$ and let $({\mathcal S}_0,{\mathcal T}_0)$ be a measurable
space such that ${\mathcal T}$ is countably generated.

\begin{theorem}{\textup{(\cite{alsmeyer}, Theorem 1)}}\label{renewala}
Let $(x_n)_{n \in {\mathbb Z}}$ be a Harris recurrent Markov chain on
$({\mathcal S}_0,{\mathcal T}_0)$ with stationary distribution $\pi$
and let $(x_n,q_n)_{n \in{\mathbb Z}}$ be an associated with it MMP on
$({\mathcal S}_0 \times{\mathbb R}, {\mathcal T}_0 \times{\mathcal B})$
such that $\mu_0:=E(q_n)>0$ and the process $q_n$ is nonarithmetic
relative to $(x_n)_{n \in{\mathbb Z}}.$ Further, let
$V_n=\sum_{i=0}^{n-1} q_i$ and let $g \dvtx  {\mathcal S}_0
\times{\mathbb R}\to{\mathbb R}$ be any measurable function satisfying
%
%e2.13 ###
\begin{equation}\label{al1}
\mbox{for } \pi\mbox{-a.e. } z \in {\mathcal S}_0, g(z,\cdot) \mbox{ is
Lebesgue-a.e. continuous},
\end{equation}
and
%
%e2.14 ###
\begin{equation}\label{al5} \int_{{\mathcal S}_0}  \sum_{n \in
{\mathbb Z}} \sup_{n \delta\leq t< (n+1) \delta} |g(z,t)|
\pi(dz)<\infty\qquad\mbox{ for some } \delta>0.
\end{equation}
Then,
%
%e2.15 ###
\begin{equation}\label{expo}
\lim_{t \to\infty} E_z^- \Biggl(\,\sum_{n=0}^\infty g (x_{n-1},t-V_n )
\Biggr)=\frac{ 1}{ \mu_0 }\int_{{\mathcal S}_0} \int_{\mathbb
R}g(u,v)\,dv\, \pi(du),
\end{equation}
for $\pi$-almost every $z \in{\mathcal S}_0.$
\end{theorem}

Under the assumptions of Theorem~\ref{renewala}, let
$\sigma_{-1}=-1,V_{-1}=0,$ and for \mbox{$n \geq0,$} let $\sigma_n=
\inf\{i>\sigma_{n-1}\dvtx V_i>V_{\sigma_{n-1}}\}$ be the ladder indexes
of the random walk~$V_n.$ Set $\widetilde V_n=V_{\sigma_n}.$ Further,
for $n \geq0$ let $\tilde x_n=x_{\sigma_n-1}$ and $\tilde
q_n=\widetilde V_n -\widetilde V_{n-1}$ ($\tilde
q_0=\sum_{i=0}^{\sigma_0-1} q_i$ and $\tilde
q_n=\sum_{i=\sigma_{n-1}}^{\sigma_n-1}q_n$ for $n \geq1$). Denote by
$\pi_1$ the unique stationary measure of the Markov chain $(\tilde
x_n)_{n \geq0}$ (existing by \cite{alsmeyer}, Theorem~2) and by $H_1$
the transition kernel of $(\tilde x_n,\tilde q_n)_{n \geq0}.$

For $t>0,$ set $\upsilon(t)=\inf\{n \geq0\dvtx  V_n>t\},$
$Z(t)=x_{\upsilon(t)-1},$ and $W(t)=V_{\upsilon(t)}-t.$ Note that
$v(t)$ is a member of the sequence $(\sigma_n)_{n \geq0}.$

\begin{corollary}[\textup{(\cite{alsmeyer}, Corollary~2)}]\label{renewal1}
Let $(x_n, V_n)_{n \geq0}$ be as in Theorem~\textup{\ref{renewala}}.
Then, with $\mu_1:= \int_{{\mathcal S}_0} E_x^- (\tilde q_0) \pi_1(d
x),$
\begin{eqnarray*}&&\lim_{t \to\infty} E_z^- (
g(Z(t),W(t)) )
\\
&& =\frac{ 1}{ \mu_1 } \int_{{\mathcal S}_0} \int
_{{\mathcal S}_0 \times(0,\infty)} \int_{[0,s)} g(v,w)\,d w\, H_1(u,dv
\times d s)\pi_1(du),
\end{eqnarray*}
holds for $\pi_1$-a.e. $z \in{\mathcal S}_0$ and for every measurable
function \mbox{$g\dvtx  {\mathcal S}_0 \times[0,\infty) \to{\mathbb
R}$} such that the function $b(z,y):=E_z^- (g (\tilde x_0,\tilde
q_0-y)\mathbf{1}_{\{\tilde q_0>y\}} ) $ satisfies \textup{(\ref
{al1})}\break and~\textup{(\ref{al5})}.
\end{corollary}

Theorem~\ref{renewala} will be applied in Section \ref{arbitar} to the
underlying Markov chain $(x_n)_{n \in {\mathbb Z}}$ restricted to the
space $({\mathcal S}_0,{\mathcal T}_0)$ defined in Section~\ref{bagr}
and to the random walk $V_n=\sum_{i=0}^{n-1} \log|\rho_i|.$ In order to
enable the application of the renewal theorem, we use a standard change
of measure argument (involving a \textit{ similarity transform} of the
transition kernel $H$) which defines a new stationary measure
$\widetilde P$ for the MMP $(x_n,\zeta_n)_{n \in{\mathbb Z}}$ under
which the Markov random walk $V_n=\sum_{i=0}^{n-1} \log|\rho_i|$ has
positive drift, that is, the expectation $\widetilde E(\log|\rho_0|)$
with respect to $\widetilde P$ is strictly positive.

We next proceed with the construction of the measure $\widetilde P.$
Observe that in virtue of Lemma~\ref{urg5}, $r_\kappa=1,$ where
$r_\kappa$ is the spectral radius of the kernel $H_\kappa$ on
$({\mathcal S}_0,{\mathcal T}_0)$ defined in \textup{(\ref{kbeta})}.
Therefore, by Proposition~\ref {urg}, there exists a positive
measurable function $h(x)\dvtx  {\mathcal S}_0 \to {\mathbb R}$ bounded
away from zero and infinity such that
%
%e2.16 ###
\begin{equation}\label{hfun}
h(x)=\int_{{\mathcal S}_0} H_\kappa(x,dy) h(y).
\end{equation}
Let $\zeta_n=(\xi_n,\rho_n),$ $n \in{\mathbb Z},$
%
%e2.17 ###
\begin{equation}\label{witih}
\widetilde H(x,dy):=\frac{ 1}{ h(x) } H_\kappa(x,dy)h(y),
\end{equation}
and let $\widetilde P$ be the stationary law of the Markov chain
$(x_n,\zeta_n)_{n \in{\mathbb Z}}$ on ${\mathcal S}\times{\mathbb R}^2$
with transition kernel
\[
\widetilde P \bigl(y_0 \in A \times B | \sigma(y_i\dvtx i<0)
\bigr)=\int_A \widetilde H(x,dz){\mathbb G}(x,z,B) |_{x=x_{-1}},
\]
where $A \in{\mathcal T}_0, B \in{\mathcal B}^{\otimes2}$ and ${\mathbb
G}(x,z,\cdot)=P(\zeta_n \in\cdot|x_{n-1}=x,x_n=z).$ That is, the law of
$(\zeta_n)_{n \in{\mathbb Z}}=(\xi_n,\rho_n)_{n \in {\mathbb Z}}$
conditioned upon $(x_n)_{n \in{\mathbb Z}}$ is the same under $P$ and
$\widetilde P,$ whereas the chain $(x_n)_{n \in{\mathbb Z}}$ has
transition kernels $H$ and $\widetilde H$, respectively. We will denote
by $\widetilde E$ the expectation with respect to $\widetilde P$ and
will use the notation
%
%e2.18 ###
\begin{equation}\label{nots}
\widetilde P^{-}_x(\cdot):=\widetilde
P(\cdot|x_{-1}=x)\quad\mbox{and}\quad\widetilde P_x(\cdot):=\widetilde
P(\cdot|x_0=x),
\end{equation}
and, correspondingly, $\widetilde E_x^{-}(\cdot):=\widetilde
E(\cdot|x_{-1}=x)$ and $ \widetilde E_x(\cdot):=\widetilde
E(\cdot|x_0=x).$

Let
%
%e2.19 ###
\begin{equation}
\label{ciah} c_h:=\sup_{x,y \in{\mathcal S}_0} h(x)/h(y).
\end{equation}
Since $c_h \in(0,\infty)$ and $c_h^{-1} H(x,A) \leq\widetilde H(x,A)
\leq c_h H(x,A),$ we have:
\begin{itemize}
\item Conditions (A1)--(A3) of Assumption~\ref{measure1} hold for the
kernel $\widetilde H.$

\item The Markov chain $(x_n)_{n \in{\mathbb Z}}$ on $({\mathcal
S}_0,{\mathcal T}_0)$ with the kernel $\widetilde H$ is Harris
recurrent and the minorization condition \textup{(\ref{lem-lem})} holds
in the following form:
%
%e2.20 ###
\begin{equation}\label{lem-lem1}
\widetilde H^m(x,A)
> r c_h^{-1} \psi(A)\qquad\forall  x \in{\mathcal S}_1,  A \in{\mathcal T}.
\end{equation}

\item The invariant measure $\pi_h$ of the kernel $\widetilde H$ is
equivalent to $\pi$ (this follows, for example, from~\cite{nummelin},
Proposition~2.4).

\item Assumptions (A7) and (\ref{calso}) hold for the sequence
$(\xi_{n},\rho_{n})_{n \in{\mathbb Z}}$ under the measure $\widetilde
P.$
\end{itemize}

\begin{lemma}\label{alemma3}
Let Assumption~\textup{\ref{measure1}} hold. Then $\widetilde E (\log|\rho_0|)>0.$
\end{lemma}

\begin{pf}
Let $V_0=0$ and
%
%e2.21 ###
\begin{equation}\label{vik} V_n=\sum_{i=0}^{n-1} \log|\rho_i|,\qquad n
\in{\mathbb N}.
\end{equation}
With $c_h$ defined in \textup{(\ref{ciah})} we obtain for any $x
\in{\mathcal S}_0$ and $\gamma>0,$
%
%e2.22 ###
\begin{eqnarray}\label{contra}
\widetilde P_x^{-} (e^{V_n} \leq e^{-\gamma n^{1/4}} )&=& \frac{ 1}{
h(x) } E_x^{-} \bigl(e^{\kappa V_n} h(x_{n-1});e^{V_n} \leq e^{-\gamma
n^{1/4}} \bigr) \nonumber
\\[-8pt]
\\[-8pt]
\nonumber &\leq& c_h E_x^{-} (e^{\kappa V_n};e^{V_n} \leq e^{-\gamma
n^{1/4}} ) \leq c_h e^{-\kappa\gamma n^{1/4}}.
\end{eqnarray}
Thus, $\lim_{n \to\infty} \widetilde P_x^{-} (V_n \leq-\gamma n^{1/4}
)=0,$ implying by the ergodic theorem that $\widetilde E(\log|\rho_0|)
\geq0.$

It remains to show that $\widetilde E(\log|\rho_0|)=0$ is impossible.
For any $ x \in{\mathcal S}_0,$ $\delta>0,$ and $\beta\in(0,\kappa)$ we
get, using Chebyshev's inequality,
\begin{eqnarray*}
\widetilde P_x^{-} (|V_n| \leq\delta n)&=&  \frac{ 1}{ h(x) } E_x^{-}
\bigl(e^{\kappa V_n} h(x_{n-1});V_n \in[-\delta n, \delta n] \bigr)
\\
&\leq& c_h e^{\kappa\delta n} P_x^{-} (V_n \geq- \delta n ) \leq c_h
e^{(\kappa+\beta) \delta n } E_x^{-} \Biggl(\, \prod_{i=0}^{n-1}
|\rho_i|^\beta\Biggr).
\end{eqnarray*}
It follows from Lemma \ref{urg5} that for all $\delta>0$ small enough
and some suitable constants $A,b>0$ that depend on $\delta,$
%
%e2.23 ###
\begin{equation}\label{nest} \sup_{ x
\in{\mathcal S}_0} \widetilde P_x^{-} (|V_n| \leq\delta n)\leq A e^{-
bn}.
\end{equation}
Therefore, the ergodic theorem implies that $\widetilde E (\log
|\rho_0|)>0.$
\end{pf}

%s3 ###
\section{Reduction to a renewal theorem}\label{proofs}
The main goal of this section is to prove the following
Proposition~\ref{gop} which reduces the limit problem for the tail of
the random variable $R$ to a renewal theorem [namely, to the checking
that \textup{(\ref{hpl})} below indeed holds a.s.]. Furthermore, some
useful estimates are obtained here and collected in Lemma
\ref{alemma1}.

Let $\Pi_0=1$ and for $n \geq1,$
%
%e3.1 ###
\begin{equation}\label{kohav}
\Pi_n=\prod_{i=0}^{n-1} \rho_i.
\end{equation}
That is $\Pi_n=\gamma_{n-1}e^{V_n}$ where $V_n$ is defined
in~\textup{(\ref{vik})} and
%
%e3.2 ###
\begin{equation}\label{gamma-s}
\gamma_n:=\operatorname{sign}(\Pi_{n+1}),\qquad n
\geq-1.
\end{equation}

\begin{proposition}\label{gop}
Let Assumption~\textup{\ref{measure1}} hold. Further, let the set
${\mathcal S}_0$ and the measure $\pi$ be as in
Lemma~\textup{\ref{urg3}} and assume that for some $\eta\in\{-1,1\}$
the following limit exists for $\pi$-almost every $z \in{\mathcal
S}_0$:
%
%e3.3 ###
\begin{equation}
\label{hpl} \widetilde K_\eta(z):=\lim_{t \to\infty} \widetilde E_z^{-}
\Biggl(\,\sum_{i=0}^\infty g_{\eta\gamma_{i-1}}(x_{i-1},t-V_i) \Biggr),
\end{equation}
where the expectation is taken according to the measure $\widetilde
P_z^{-}$ defined in \textup{(\ref{nots})} and the nonnegative functions
$g_\gamma\dvtx  {\mathcal S}_0 \times{\mathbb R}\to[0,\infty)$ are
defined for $\gamma \in\{-1,1\}$ by
%
%e3.4 ###
\begin{equation}\label{gi1}
g_\gamma(x,t)=\frac{e^{-t}}{h(x)}\int_0^{e^t} v^\kappa \bigl[P_x^{-}
(\gamma R>v )- P_x^{-} \bigl(\gamma(R-\xi_0)>v \bigr) \bigr]\,dv.
\end{equation}
Then, for $\pi$-almost every $z \in{\mathcal S}_0, \lim_{t \to \infty}
t^\kappa P_z^{-}(\eta R>t)=h(z) \widetilde K_\eta(z).$
\end{proposition}

We note that certain particular cases of this proposition are the basis
for the proofs in \cite{saporta} and in \cite{mars}. All these results
are adaptations to various Markovian situations of a particular case of
the ``implicit renewal'' theorem of Goldie (cf. Theorem 2.3 in
\cite{goldie}). For the sake of completeness, a proof of the
proposition is provided at the end of this section.

We begin by proving the following technical lemma:

\begin{lemma}\label{alemma1}
Let Assumption~\textup{\ref{measure1}} hold. Then the following
assertions hold true:
\begin{longlist}[(a)]
\item[(a)] There exists constants $M_g>0$ and $\varepsilon_g>0$ such
that for $\pi$-almost every $x \in{\mathcal S}_0,$
%
%e3.5 ###
\begin{equation}\label{la1} |g_\eta(x,t)| \leq M_g
e^{-\varepsilon_g |t|},
\end{equation}
for any $t \in{\mathbb R}$ and $\eta\in \{-1,1\}.$

In particular, for any $\delta>0$ there exists a constant $M(\delta)>0$
such that
%
%e3.6 ###
\begin{equation}\label{del}
\sum_{n \in{\mathbb Z}}   \sup_{n \delta\leq t< (n+1) \delta} \biggl\{
\max_{\eta\in\{-1,1\}} |g_\eta(x,t)|  \biggr\} \leq M(\delta)
\end{equation}
for $\pi$-almost every $x \in{\mathcal S}_0.$

\item[(b)] For any $\delta>0$ there exists a constant
$M_u=M_u(\delta)>0$ such that
\begin{eqnarray*}\sum_{i=0}^\infty\sup_{z \in{\mathcal S}_0}
\widetilde P_z^{-} (V_i \in [-\delta,\delta]) \leq M_u.
\end{eqnarray*}

\item[(c)] There exists a constant $M_r>0$ such that, for $\pi$-almost
every $z \in{\mathcal S}_0,$
%
%e3.7 ###
\begin{equation}
\label{fubini} \sum_{i=0}^\infty\widetilde E_z^{-} \Bigl( \max_{\eta\in
\{-1,1\}} |g_{\eta}(x_{i-1},t-V_i) |   \Bigr) \leq M_r\qquad \forall  t \in
{\mathbb R}.
\end{equation}
\end{longlist}
\end{lemma}

\begin{pf}
(a) First, assume that $t>0.$ Let
%
%e3.8 ###
\begin{equation}\label{ci-tilde}
\tilde c_h=\max_{x \in S_0} 1/h(x).
\end{equation}
For any $\varepsilon\in(0,1),$ we get from \textup{(\ref{gi1})}:
%
%e3.9 ###
\begin{eqnarray}\label{above}
\nonumber |g_\eta(x,t)| &\leq& \tilde c_h e^{-t}\int_0^{e^t}
v^\kappa\big| P_x^{-} (\eta R>v)-P_x^{-} \bigl(\eta (R-\xi_0)>v\bigr)
\big| \,dv
\\
&\leq& \tilde c_h e^{-\varepsilon t}\int_0^{e^t} v^{\kappa
+\varepsilon-1} \big| P_x^{-} (\eta R>v)- P_x^{-}
\bigl(\eta(R-\xi_0)>v\bigr) \big| \,dv
\\
\nonumber &\leq& \tilde c_h \kappa^{-1} e^{-\varepsilon t} E_x^{-}
\bigl( | [(\eta R)^+]^{\kappa+\varepsilon}- [(\eta R- \eta
\xi_0)^+]^{\kappa+\varepsilon} |  \bigr),
\end{eqnarray}
where the last inequality follows from \cite{goldie}, Lemma 9.4.

To bound the right-hand side in \textup{(\ref{above})} we will exploit
an argument similar to the proof of \cite{goldie}, Theorem~4.1. We
have,
\begin{eqnarray*}
|g_\eta(x,t)| \leq\tilde c_h \kappa^{-1} e^{-\varepsilon t}
[I_1(x)+I_2(x)+I_3(x)+I_4(x)],
\end{eqnarray*}
where
\begin{eqnarray*}
I_1(x)&:=& E_x^{-} \bigl( \mathbf{1}_{\eta\xi_0 <\eta R \leq0 } (\eta
R-\eta \xi_0)^{\kappa+\varepsilon} \bigr),
\\
I_2(x)&:=& E_x^{-} \bigl( \mathbf{1}_{0<\eta R \leq\eta\xi_0 } (\eta
R)^{\kappa+\varepsilon } \bigr),
\\
I_3(x)&:=& E_x^{-} \bigl( \mathbf{1}_{\eta R>0, \eta\xi_0<0} [(\eta
R-\eta \xi_0)^{\kappa+\varepsilon} -(\eta R)^{\kappa+\varepsilon} ]
\bigr),
\\
I_4(x)&:=& E_x^{-} \bigl( \mathbf{1}_{0 \leq\eta\xi_0< \eta R} [(\eta
R)^{\kappa+\varepsilon} -(\eta R-\eta\xi_0)^{\kappa+\varepsilon} ]
\bigr).
\end{eqnarray*}
It follows from \textup{(\ref{calso})} that the sum $I_1(x)+I_2(x)$ is
bounded by $c_\xi^{\kappa+\varepsilon}.$ It remains therefore to bound
$I_3(x)$ and $I_4(x).$ For this purpose we will use the following
inequalities valid for any $\gamma>0$ and $A>0,B>0$ (this is exactly
(9.26) and (9.27) in \cite{goldie}):
\[
(A+B)^\gamma\leq2^\gamma(A^\gamma+B^\gamma)
\]
and
\[
(A+B)^\gamma- A^\gamma\leq \cases{ B^\gamma,& \quad if $0
\leq\gamma\leq1$, \cr \gamma B (A+B)^{\gamma-1}, &\quad if $\gamma>1$.}
\]
We obtain that $I_3(x)+I_4(x)\leq a_\varepsilon,$ where
\[
a_\varepsilon:= \cases{ c_\xi^{\kappa+\varepsilon}, &\quad if  $\kappa
+\varepsilon\leq1$, \cr (\kappa+\varepsilon)
c_\xi2^{\kappa+\varepsilon-1}E_x (|R|^{\kappa +\varepsilon
-1}+c_\xi^{\kappa+\varepsilon-1} ), &\quad if $\kappa+\varepsilon>1$.}
\]
By Lemma~\ref{urg5} and the ellipticity condition
\textup{(\ref{calso})}, for any $\delta>0$ small enough, there exists a
constant $L_\delta>0$ independent of $x$ such that
%
%e3.10 ###
\begin{equation}\label{lde} E_x^{-} (
|R|^{\kappa-\delta} ) \leq L_\delta\qquad \forall x \in{\mathcal S}_0.
\end{equation}
This yields \textup{(\ref{la1})} for all $t>0$ and appropriate
constants $M_g,\varepsilon_g>0$ that do not depend on $t.$

Further, \textup{(\ref{gi1})} implies that $|g_\eta(x,0)| \leq \tilde
c_h,$ where the constant $\tilde c_h$ is defined in~\textup{(\ref
{ci-tilde})}, and that for $t<0$ and any $\varepsilon\in(0,\kappa),$
\begin{eqnarray*}
|g_\eta(x,t)| &\leq& \tilde c_h e^{-t}\int_0^{e^t} v^\kappa\big|
P_x^{-} (\eta R>v)-P_x^{-} \bigl(\eta(R-\xi_0)>v\bigr) \big| \,dv
\\
&\leq& \tilde c_h e^{\varepsilon t}\int_0^{e^t} v^{\kappa
-\varepsilon-1} \big| P_x^{-} (\eta R>v)-P_x^{-} \bigl(\eta(R-\xi_0)>v\bigr) \big|\, dv
\\
&\leq& \tilde c_h \kappa^{-1} e^{\varepsilon t} E_x^{-} \bigl(  |
[(\eta R)^+]^{\kappa-\varepsilon}- [(\eta R- \eta\xi_0)^+]^{\kappa
-\varepsilon} | \bigr),
\end{eqnarray*}
where the last inequality follows, similarly to \textup{(\ref{above})},
from \cite{goldie}, Lemma~9.4. Thus, for $t<0,$
\[
|g_\eta(x,t)| \leq\tilde c_h \kappa^{-1} e^{-\varepsilon|t|} E_x^{-}
\bigl( |R| ^{\kappa-\varepsilon}+ (|R|+c_\xi)^{\kappa-\varepsilon}
\bigr).
\]
This completes the proof in view of \textup{(\ref{lde})}.

(b) Follows from \textup{(\ref{nest})}, since $\widetilde P_x^{-} ( V_i
\in[-\delta, \delta] ) \leq \widetilde P_x^{-} ( V_i \in[-i\delta,
i\delta] )$ for any $x \in{\mathcal S}_0$ and $i \in{\mathbb N}.$

(c) Fix any $\delta>0$ and denote for $t \in{\mathbb R}$ and $n
\in{\mathbb Z},$ $I_{n,\delta}^t=[t+n \delta , t+(n+1) \delta).$ Then,
it follows from the previous parts of the lemma that
\begin{eqnarray*}
&& \sum_{i=0}^\infty\widetilde E_z^{-} \bigl( |g_\eta (x_{i-1},t-V_i) |
\bigr)
\\
&&\qquad =\sum_{i=0}^\infty\int_{{\mathcal S}_0} \int_{\mathbb R}
|g_\eta(x,t-s) | \widetilde P_z^{-} (x_{i-1} \in dx,V_i \in d s)
\\
&&\qquad \leq\sum_{n \in{\mathbb Z}} \sup_{x \in{\mathcal S}_0 , s \in
I_{n,\delta}^t} |g_\eta(x,t-s) | \sum_{i=0}^\infty
 \sup_{z \in{\mathcal S}_0} \widetilde P_z^{-} (V_i \in t-I_{n,\delta
}^s)
\\
&&\qquad \leq M(\delta) \sum_{i=0}^\infty  \sup_{z \in{\mathcal S}_0}
\widetilde P_z^{-} (V_i \in[-\delta,\delta]) \leq M(\delta) \cdot M_u,
\end{eqnarray*}
where the last but one inequality follows from \textup{(\ref{del})} and
from the fact that $\sup_{z \in{\mathcal S}_0}
\sum_{i=0}^\infty\widetilde P_z^{-} (V_i \in t-I_{n,\delta}^s)
\leq\sup_{z \in{\mathcal S}_0} \sum_{i=0}^\infty \widetilde P_z^{-}
(V_i \in[-\delta,\delta])$ (cf. \cite{alsmeyer},\break Lemma~A.2).
\end{pf}

\begin{pf*}{Proof of Proposition~\ref{gop}}
Let $U_0=R,$ and for $n \geq1,$
%
%e3.11 ###
\begin{equation}\label{kohav11}
R_n=\sum_{i=0}^{n-1} \xi_i \Pi_i,\qquad  U_n=(R-R_n)/\Pi_n,
\end{equation}
where $\Pi_n$ are defined in \textup{(\ref{kohav})}. Recall that
$\Pi_n=\gamma_{n-1}e^{V_n}.$ Following Goldie \cite{goldie}, we write
for any numbers $n \in{\mathbb N},$ $t \in{\mathbb R},$
$\eta\in\{-1,1\}$ and any $z \in{\mathcal S}_0,$
\begin{eqnarray*}
P_z^{-}(\eta R>e^t)&=&\sum_{i=0}^{n-1} [
P_z^{-}(\eta\gamma_{i-1}e^{V_i}
U_i>e^t)-P_z^{-}(\eta\gamma_ie^{V_{i+1}} U_{i+1}>e^t)]
\\
&&{} + P_z^{-}(\eta \Pi_n U_n>e^t),
\end{eqnarray*}
where the random variable $U_i$ is defined in \textup{(\ref{kohav11})}.

For $n \geq-1$ let $\hat x_n=(x_n,\gamma_n)$ and $\Omega={\mathcal
S}\times\{-1,1\}\times{\mathbb R}.$ To shorten the notation, we denote
$\sum_{\gamma\in\{-1,1\}}\int _{\mathcal S}\int_{\mathbb
R}F(\gamma,x,u) \mu(\gamma_{i-1}=\gamma,x_i \in dx,V_i \in du)$ by
$\int_\Omega F(\gamma,x,u) \mu(\hat x_i \in(dx,\gamma),V_i \in du)$ for
a measurable function $F$ and a probability measure $\mu$ on $\Omega.$
We have, using the identity $U_i=\xi_i+\rho_i U_{i+1},$
\begin{eqnarray*}
&& P_z^{-}(\eta\gamma_{i-1}e^{V_i} U_i>e^t)-P_z(\eta
\gamma_ie^{V_{i+1}} U_{i+1}>e^t)
\\
&&\qquad = \int_\Omega P\bigl(\eta\gamma U_i>e^{t-u}|\hat
x_{i-1}=(x,\gamma) , V_i=u\bigr)P_z^{-}\bigl(\hat x_{i-1}
\in(dx,\gamma),V_i\in du\bigr)
\\
&&\qquad\quad{} - \int_\Omega
P\bigl(\eta\gamma\rho_iU_{i+1}>e^{t-u}|\hat x_{i-1}=(x,\gamma),
V_i=u\bigr)
\\
&&\phantom{\qquad\quad{} - \int_\Omega} {}\times P_z^{-}\bigl(\hat
x_{i-1}\in(dx,\gamma), V_i\in du\bigr)
\\
&&\qquad =\int_\Omega e^{-\kappa(t-u)}f_{\eta\gamma}(x,t-u)
P_z^{-}\bigl(\hat x_{i-1} \in(dx, \gamma),V_i\in du\bigr),
\end{eqnarray*}
where we denote
\[
f_\gamma(x,t)=e^{\kappa t}\bigl[P_x^{-} (\gamma
R>e^t)-P_x^{-} \bigl(\gamma (R-\xi_0)>e^t\bigr)\bigr]\qquad\mbox{for }
\gamma\in\{-1,1\}.
\]
Thus, letting $\delta_n(z,\eta,t)= e^{\kappa t} P_z^{-}(\eta
\gamma_{n-1}e^{V_n} U_n>e^t)$ we obtain
\begin{eqnarray*}
\check{r} _z(\eta,t)&:=&e^{\kappa t}P_z^{-}(\eta R
>e^t)
\\
&\hspace*{3pt}  =& \sum_{i=0}^\infty\int_\Omega f_{\eta\gamma} (x,t-u)
e^{\kappa u} P_z^{-}\bigl(\hat x_{i-1} \in(dx,\gamma), V_i\in du\bigr)
+ \delta_n(z,\eta,t)
\\
&\hspace*{3pt}  =& \sum_{i=0}^{n-1} \int_\Omega f_{\eta\gamma} (x,t-u)
\frac{h(z)}{h(x)} \widetilde P_z^{-}\bigl(x_{i-1} \in(dx,\gamma),
V_i\in du\bigr)+ \delta_n(z,\eta,t).
\end{eqnarray*}
We have $P (\lim_{n\to\infty} \delta_n(z,\eta,t)=0 )=1$ for any fixed
$t>0,$ $\eta\in\{-1,1\},$ and $z \in{\mathcal S}_0,$ because $P$-a.s.,
$\Pi_n U_n \to0$ as $n$ goes to infinity. Therefore $P$-a.s.,
\begin{eqnarray*}r_z(\eta,t)=\sum_{i=0}^\infty
\int_\Omega f_{\eta\gamma}(x,t-u) \frac{h(z)} {h(x)} \widetilde
P_z^{-}\bigl(\hat x_{i-1} \in(dx,\gamma),V_i\in du\bigr).
\end{eqnarray*}
We will use the following Tauberian lemma:

\begin{lemma}[\textup{(\cite{goldie}, Lemma~9.3)}]\label{igol}
Let $R$ be a random variable such that for some constants $\kappa>0$
and $K \geq0,$ $\lim_{t \to\infty}t^{-1} \int_0^t u^\kappa P(R>u)\,d
u=K.$ Then $\lim_{t \to\infty}t^\kappa P(R>t)=K.$
\end{lemma}

It follows from Lemma~\ref{igol} that in order to prove that for some
$\eta\in\{-1,1\},$ the limit $\lim_{t \to\infty}t^\kappa P(\eta R>t)$
exists and is strictly positive, it suffices to show that for
$\pi\mbox{-a.s.}$ every $z \in{\mathcal S}_0,$ there exists
%
%e3.12 ###
\begin{equation}\label{check}
\lim_{t \to\infty} \check{r}_z(\eta,t) \in(0,\infty),
\end{equation}
where the smoothing transform $\check{q}$ is defined for a measurable
function $q\dvtx  {\mathbb R}\to{\mathbb R}$ bounded on $(-\infty,t]$
for all $t$ by $\check{q}(t):=\int_{-\infty}^t e^{-(t-u)}q(u)\,du.$

For $\gamma\in\{-1,1\}$ let
\begin{eqnarray*}
g_\gamma(x,t)&:=&\frac{1}{h(x)} \int_{-\infty}^t e^{-(t-u)}
f_\gamma(x,u)\,du
\\
&\hspace*{3pt} =& \frac{1}{h(x)} \int_{-\infty}^t e^{-(t-u)} e^{\kappa
u}\bigl[P_x^{-} (\gamma R>e^u)-P_x^{-}
\bigl(\gamma(R-\xi_0)>e^u\bigr)\bigr]\,du
\\
&\hspace*{3pt} =& \frac{e^{-t}}{h(x)}\int_0^{e^t} v^\kappa\bigl[P_x^{-}
(\gamma R>v)- P_x^{-} \bigl(\gamma(R-\xi_0)>v\bigr)\bigr]\,dv.
\end{eqnarray*}
Then, using \textup{(\ref{fubini})} and the Fubini theorem, we obtain
for any $z \in{\mathcal S}_0,$
\begin{eqnarray*}
\check{r}_z(\eta,t)&=& \int_{-\infty}^t e^{-(t-w)}r_z(\eta,w)\,dw
\\
&=&\int_{-\infty}^t \!e^{-(t-w)} \sum_{i=0}^\infty\int_\Omega\!
f_{\eta\gamma}(x,w-u) \frac{h(z)} {h(x)} \widetilde P_z^{-}\bigl(\hat
x_{i-1} \in(dx,\gamma), V_i\in du\bigr)\,dw
\\
&=& \sum_{i=0}^\infty\int_\Omega g_{\eta\gamma}(x,t-u) h(z) \widetilde
P_z^{-}\bigl(\hat x_{i-1} \in(dx,\gamma), V_i\in du\bigr)
\\
&=& h(z) \widetilde E_z^{-} \Biggl(\,  \sum_{i=0}^\infty g_{\eta
\gamma_{i-1}}(x_{i-1},t-V_i)   \Biggr).
\end{eqnarray*}
This completes the proof of Proposition~\ref{gop}.
\end{pf*}

%s4 ###
\section{The auxiliary Markov chain $\hat x_n=(x_n,\gamma_n)$}
\label{aux} To deal with the case where $P(\rho_0<0)>0$ we introduce
the Markov chain $\hat x_n=(x_n,\gamma_n),$ $n \geq-1,$ where the
random variables $\gamma_n$ are defined in \textup{(\ref{gamma-s})}. It
will turn out (cf. Proposition~\ref{still}) that the space ${\mathcal
S}_0 \times\{-1,1\}$ can be partitioned into at most two measurable
subsets such that the restriction of $\hat x_n$ to either one of them
satisfies Assumption~\ref{measure1}. Therefore, the Markov renewal
theorem (Theorem~\ref{renewala}) can be applied to the
irreducible\vadjust{\goodbreak} components of the MMP $(\hat
x_n,\log|\rho_n|).$ This fact is the key to the proof (given in the
next section) that the limit in \textup{(\ref{hpl})} exists $\pi$-a.s.
and has the properties stated in Theorem~\ref{main-rec-3}.

Let $\widehat H$ be the transition kernel of $\hat x_n$ on the product
space ${\mathfrak S}:={\mathcal S}_0 \times\{ -1,1\},$ and let
$\hat\pi$ be the probability measure on ${\mathfrak S}$ defined by
$\hat\pi(A \times\eta)=1/2  \pi(A)$ for any $\eta\in\{-1,1\}$ and $A
\in{\mathcal T}_0.$ It is easy to see that $\hat\pi$ is a stationary
distribution of the Markov chain $\hat x_n.$

\begin{proposition}
\label{still} Let Assumption~\textup{\ref{measure1}} hold and suppose
in addition that $P(\rho_0<0)>0.$ Then, there exist two disjoint
measurable subsets ${\mathfrak S}_1$ and ${\mathfrak S}_{-1}$ of
${\mathfrak S}$ such that:
\begin{longlist}[(iii)]
\item[(i)] Either $\hat\pi({\mathfrak S}_1)=\hat\pi({\mathfrak
S}_{-1})=1/2,$ or ${\mathfrak S}_1=\varnothing$ and ${\mathfrak
S}_{-1}={\mathfrak S}.$

\item[(ii)] $\widehat H(\hat x,{\mathfrak S}_n)=1$ for every $\hat x
\in{\mathfrak S}_n,$ $n=-1,1.$

\item[(iii)] ${\mathfrak S}_1=\varnothing$ if and only if Condition \textup{G}
is satisfied.

\item[(iv)] \textup{(A1)--(A3)} of Assumption~\textup{\ref{measure1}}
hold for the Markov chain $(\hat x_n)_{n \geq-1}$ restricted to either
${\mathfrak S}_1$ (provided that it is not the empty set) or
${\mathfrak S}_{-1}.$
\end{longlist}
\end{proposition}

\begin{pf}
(i)--(ii) Say that for $ \hat x \in{\mathfrak S},A \in{\mathcal
T}_0,\gamma\in \{-1,1\},$
\[
\hat x \not\succ A \times \{\gamma\}\qquad\mbox{if }
\sum_{n=1}^\infty\widehat H^n(\hat x,A \times\{\gamma\})=0,
\]
and $\hat x \succ A \times\{\gamma\}$ otherwise.

Since the Markov chain $(x_n)_{n \in{\mathbb Z}}$ is $\pi
$-irreducible, for any $\hat x \in{\mathfrak S}$ and $A \in{\mathcal
T}_0$ such that $\pi(A)>0$ either $\hat x \succ A \times\{1\}$ or $\hat
x \succ A \times\{-1\}.$ For $\hat x \in{\mathfrak S}$ and
$\eta\in\{-1,1\}$ let:
\[
\digamma_\eta(\hat x)= \bigl\{A \in{\mathcal T}_0\dvtx  \pi(A)>0
\mbox{ and } \hat x \not\succ A \times \{\eta\}\bigr\},
\]
and set $\digamma(\hat x)=\digamma_1(\hat x) \cup \digamma_{-1} (\hat
x).$ Note that $\digamma_1(\hat x) \cap \digamma_{-1}(\hat
x)=\varnothing$.

Roughly speaking, the set ${\mathfrak S}_1$ is defined below as an
element of $\digamma(x^*)$ of maximal $\hat\pi$-measure for some $x^*
\in {\mathfrak S},$ and ${\mathfrak S}_{-1}$ as its complement in
${\mathfrak S}.$

To be precise, let
\[
\varsigma_\eta(\hat x)=\sup\{\pi(A)\dvtx  A \in \digamma_\eta(\hat x)
\},\qquad \eta\in\{-1,1\},\hat x \in{\mathfrak S},
\]
and $\varsigma(\hat x)=\varsigma_{-1}(\hat x)+\varsigma_1(\hat x).$ If
$\varsigma(\hat x)=0$ for every $\hat x \in{\mathfrak S},$ set
${\mathfrak S}_1=\varnothing$ and \mbox{${\mathfrak S}_{-1}={\mathfrak
S}.$} Conclusions (i)--(ii) follow trivially in this case, in
particular the chain $(x_n,\gamma_n)$ is $\hat\pi$-irreducible.

Assume now that $\varsigma(x^*)>0$ for some $x^* \in{\mathfrak S}.$ We
will next construct two sets~$A_\eta,$ $\eta\in\{-1,1\},$ such that
$A_\eta\in\digamma_\eta(x^*)$ and $\pi(A_\eta)= \varsigma_\eta(x^*).$
We will then show that $\varsigma(x^*)=\pi(A_1)+\pi(A_{-1})=1$ and will
define (up to a $\hat\pi$-null set) ${\mathfrak S}_1:=  (A_{-1} \times
\{-1\} ) \cup( A_1 \times\{1\} ).$

For $\eta\in\{-1,1\},$ let $A_{\eta,n} \in\digamma_\eta(x^*),$ $n
\in{\mathbb N},$ be a sequence of [empty if
\mbox{$\varsigma_\eta(x^*)=0$}] sets in $\digamma_\eta(x^*)$ such that
$\pi(A_{\eta,n})>\varsigma_\eta(x^*)-1/n$ for any $n \in{\mathbb N},$
and define $A_\eta=\bigcup_{n=1}^\infty A_{\eta,n}.$ Since the
collections of sets $\digamma_\eta(x^*)$ are closed
with\vadjust{\goodbreak} respect to countable unions,
$A_\eta\in\digamma_\eta(x^*)$ and
$\pi(A_1)+\pi(A_{-1})=\varsigma(x^*).$ Put $A_0=A_{-1} \cup A_1,$
$B_0=S_0-A_0,$ and set
\begin{eqnarray*}
\widetilde{\mathfrak S}_1 &=&  (A_{-1} \times\{-1\} ) \cup( A_1
\times\{1\} ),
\\
\widetilde{\mathfrak S}_{-1}&=&  (A_{-1} \times\{1\} ) \cup(A_1
\times\{-1\} ) \cup( B_0\times\{1\} ) \cup(B_0 \times\{-1\} ).
\end{eqnarray*}
Thus, $\widetilde{\mathfrak S}_{-1}$ is the complement of
$\widetilde{\mathfrak S}_1$ in the set ${\mathfrak S}={\mathcal S}_0
\times\{-1,1\}.$ Since $A_\eta\in\digamma_\eta(x^*)$ is the maximal set
such that $x^* \not\succ A_\eta\times\{\eta\},$ it follows immediately
that $x^* \not\succ A \times\{\eta\}$ and $x^* \succ A \times\{-\eta\}$
for any $\pi$-positive $A \subset A_\eta.$

We will now show, using the irreducibility of the Markov chain
$(x_n)_{n \in{\mathbb Z}},$ that
%
%e4.1 ###
\begin{equation}\label{from} \hat\pi(N_{-1}
\cup N_1 )=0\qquad\mbox{where }  N_\eta:=\bigl\{\hat x \in\widetilde
{\mathfrak S}_{-1}\dvtx  \hat x \succ A_\eta\times\{\eta\}\bigr\}.
\end{equation}
Note that \textup{(\ref{from})} yields $\pi(B_0)=0$ because for all $x
\in{\mathcal S}_0$ either $(x,1) \in N_\eta$ or $(x,-1) \in N_\eta,$
$\eta\in\{-1,1\}.$

To see that \textup{(\ref{from})} is true, observe that
$N_\eta=\bigcup_{m \in{\mathbb N}} \{\hat x\dvtx  \widehat H^m(\hat x,
A_\eta\times\{\eta\})>0 \}$ are measurable sets, and $\hat\pi
(N_\eta)>0$ implies that there exist $m \in{\mathbb N},$ $N_0 \in
{\mathcal T}_0,$ and $\gamma\in\{-1,1\}$ such that
%
%e4.2 ###
\begin{equation}\label{bnul}
\widehat H^m(x^*, N_0 \times\{\gamma\})>0\quad\mbox{and}\quad (x,
\gamma) \succ A_\eta \times\{\eta\}\qquad\forall  x \in N_0.
\end{equation}
But \textup{(\ref{bnul})} yields
\begin{eqnarray*}
&& \sum_{n=0}^\infty\widehat H^{m+n} (x^*,A_\eta\times\{\eta\} )
\\
&&\qquad \geq\sum_{n=0}^\infty\int_{N_0} \widehat H^m (x^*, dy \times
\{\gamma\} ) \widehat H^n \bigl((y,\gamma), A_\eta\times\{\eta\}
\bigr)>0,
\end{eqnarray*}
which is impossible since $x^* \not\succ A_\eta\times \{\eta\}$ by our
construction.

Finally, we observe that \textup{(\ref{from})} implies that
%
%e4.3 ###
\begin{equation}\label{frok}
\hat\pi( \widebar N_{-1} \cup\widebar N_1 )=0\qquad\mbox{where }
\widebar N_\eta:=\bigl\{\hat x \in\widetilde{\mathfrak S}_1\dvtx  \hat
x \succ A_\eta \times \{-\eta\}\bigr\}.
\end{equation}
Indeed, if $(x,\gamma) \in\widebar N_\eta$ then $(x,-\gamma) \in
N_\eta$ and hence $\hat\pi(\widebar N_\eta)=\hat \pi (N_\eta)=0$ for
$\eta\in\{-1,1\}.$

To complete the proof, we set
\[
{\mathfrak S}_1= (A_{-1} \times\{-1\} ) \cup( A_1 \times\{1\} )
-\widebar N_{-1} \cup\widebar N_1,
\]
and
\[
{\mathfrak S}_{-1}= (A_{-1} \times\{1\} ) \cup (A_1 \times\{-1\}
)-N_{-1} \cup N_1.
\]
Since $\pi(B_0)=0,$ \textup{(\ref{from})} and \textup{(\ref{frok})}
imply that $\hat\pi({\mathfrak S}_1)=\hat\pi ({\mathfrak S}_{-1})=1/2$
[recall that $\pi(A_1 \cap A_{-1})=0$] and that conclusion (ii) of the
proposition holds as well.

(iii) The claim is immediate from the definition of the sets $A_1$ and
$A_{-1}.$

(iv) Let $\hat\mu$ be the probability measure on ${\mathfrak S}$
defined by $\hat\mu(A \times\eta)=1/2 \mu(A),$ where $\mu(\cdot)$ is
given by assumption~(A3). Since $\widehat H^{m_1} ( (x,\gamma), A
\times\{\eta\} ) \leq\break H^{m_1}(x,A),$ it follows from~(A3) that
there exists a measurable density kernel $\hat h(\hat x, \hat y)\dvtx
{\mathfrak S}^2 \to[0,\infty)$ such that or any $\hat x \in{\mathfrak
S}, \eta\{-1,1\},A \in{\mathcal T}_0,$
%
%e4.4 ###
\begin{equation}\label{a3ho} \widehat
H^{m_1} ( \hat x,A \times\{\eta\})=\int_{A \times\{\eta\}} \hat h(\hat
x, \hat y) \hat\mu(d \hat y),
\end{equation}
and the family of functions $\{\hat h(\hat x,\cdot)\dvtx {\mathfrak
S}\to[0,\infty)\}_{\hat x \in {\mathfrak S}}$ is uniformly integrable
with respect to the measure $\hat\mu.$ Thus assumptions (A1)~and~(A3)
hold for the Markov chain $(x_n,\gamma_n)_{n \geq-1}.$ Moreover, the
Markov chain $(x_n,\gamma_n)_{n \geq-1},$ when restricted to either
${\mathfrak S}_1$ or ${\mathfrak S}_{-1},$ is clearly
$\hat\pi$-irreducible which in combination with~(\ref{a3ho})
shows~(iv).
\end{pf}

%s5 ###
\section{Distribution tail of $R$}\label{arbitar}
In this section we complete the proof of
Theorems~\textup{\protect\ref{main-rec}},
\textup{\protect\ref{main-rec-3}} and \textup{\protect\ref{mishne}}.

%s5.1 ###
\subsection[Proofs of Theorems 1.5
and 1.6 for $P(\rho_0>0)=1$]{Proofs of Theorems \textup{\protect\ref{main-rec}}
and \textup{\protect\ref{main-rec-3}} for $P(\rho_0>0)=1$} \label{mrec}
In view of Proposition~\ref{gop}, the following lemma completes the
proof of Theorem~\ref{main-rec} and of Theorem~\ref{main-rec-3} in the
case where $P(\rho_0>0)=1.$
\begin{lemma} \label{comph} Let Assumption~\textup{\ref{measure1}} hold and suppose
that $P(\rho_0>0)=1.$ Then the following assertions hold true for $\eta
\in\{-1,1\}$:
\begin{longlist}[(a)]
\item[(a)] The limit in \textup{(\ref{hpl})} exists for $\pi$-a.e. $z
\in{\mathcal S}_0$ and does not depend on $z.$

\item[(b)] If in addition $P(\xi_0>0)=1,$ then the limit is $\pi$-a.s.
strictly positive.

\item[(c)] $\pi(K_\eta(x)>0 ) \in\{0,1\}.$
\end{longlist}
\end{lemma}

\begin{pf}
(a) In view of Lemma~\ref{alemma3}, estimate \textup{(\ref{del})}, and
the properties of the measure $\widetilde P$ listed right before the
statement of Lemma~\ref{alemma3}, we can apply Theorem~\ref{renewala}
to the restriction of the underlying Markov chain $(x_n)_{n \in{\mathbb
Z}}$ on $({\mathcal S}_0,{\mathcal T}_0)$ with transition kernel
$\widetilde H,$ the associated with it random walk
$V_n=\sum_{i=0}^{n-1} \log|\rho_i|,$ and the functions $g_\eta$ defined
in \textup{(\ref{gi1})}. It follows from \textup{(\ref {expo})} that
the limit in \textup{(\ref{hpl})} is $\pi_h$-a.s. (and thus also
$\pi$-a.s.) equal to
%
%e5.1 ###
\begin{equation}\label{lhsk}
\widetilde K_\eta= \frac{ 1}{ \tilde a } \int_{{\mathcal S}_0}
\int_{\mathbb R}g_\eta(x,t)  \pi_h(dx)\,dt,
\end{equation}
where $\tilde a=\widetilde E (\log\rho_0).$

(b) It follows from Proposition~\ref{gop} and \textup{(\ref{lhsk})}
that for $\pi_h$-almost every $z \in{\mathcal S}_0$ (compare with the
formula (4.3) in \cite{goldie}),
%
%e5.2 ###
\begin{eqnarray}\label{expl}
\nonumber \hspace*{10mm} && \lim_{t \to\infty}t^\kappa
P_z^-(R>t)
\\
\nonumber &&\qquad =  h(z)\widetilde K_1(z)
\\
&&\qquad = \frac{ h(z)}{ \tilde a } \int _{{\mathcal S}_0}
\int_{\mathbb R}g_1(x,t) \pi_h(dx)\,dt
\\
\nonumber &&\qquad = \frac{ h(z)}{ \tilde a }\int_{\mathcal
S}\frac{1}{h(x)} \int_{\mathbb R}e^{-t}\int_0^{e^t} v^\kappa[P_x^{-}
(R>v)
\\
\nonumber &&\hspace*{54.5mm} {} - P_x^{-} (R-\xi_0>v)]\,dv \,dt\,
\pi_h(dx)
\\
\nonumber &&\qquad =\frac{ h(z)}{ \tilde a } \int_{\mathcal
S}\frac{1}{h(x)} \int _0^\infty v^{\kappa-1}[P_x^{-} (R>v)- P_x^{-}
(R-\xi_0>v)]\,dv\, \pi_h(dx)
\\
\nonumber &&\qquad =\frac{ h(z)}{ \tilde a \kappa} \int_{\mathcal
S}\frac{1}{h(x)} E_x^{-} [R^\kappa-(R-\xi_0)^\kappa] \pi_h(dx)>0,
\end{eqnarray}
where the last but one equality is obtained by change of the order of
the integration between $d t$ and $d v$ while the last one follows from
\cite{goldie}, Lemma~9.4. Since $\pi_h$ is equivalent to $\pi$ and
$P(R>R-\xi_0>0)=1,$ this completes the proof of the claim.

(c) The claim follows from Proposition~\ref{gop} and the fact that the
limit $\widetilde K_1$ in \textup{(\ref{hpl})} does not depend on~$z.$
\end{pf}

%
%s5.2 ###
\subsection[Proof of Theorem 1.6 for $P(\rho_0<0)>0$]{Proof of Theorem \textup{\protect\ref{main-rec-3}} for $P(\rho_0<0)>0$}\label{zvonok}
(a) Just as in the case $P(\rho_0>0),$ it follows from
Theorem~\ref{renewala}, applied separately to the irreducible
components of the Markov chain $(\hat x_n)_{n \geq-1},$ the random walk
$V_n,$ and the function $g_{\eta\gamma_n}(x_{n-1},t-V_n)$ defined in
\textup{(\ref{gi1})}, that the limits in \textup{(\ref{hpl})} and hence
in \textup{(\ref{gvul1})} exist for $\pi$-almost every $x \in{\mathcal
S}_0.$

(b)--(c) We shall continue to use the notation introduced in
Section~\ref{aux}. Similarly to \textup{(\ref{kbeta})}, define the
kernel $\widehat H_\beta(x,\cdot)$ on ${\mathfrak S}$ by
\[
\widehat H_\beta(\hat x,d \hat y)=\widehat H(\hat x,d \hat y) E (
|\rho_0|^\beta|\hat x_{-1}=\hat x,\hat x_0=\hat y ),
\]
and the function $\hat h\dvtx {\mathfrak S} \to(0,\infty)$ by the
following rule:
\[
\hat h(\hat x)=h(x)\qquad\mbox{for } x=(x,\gamma),
\]
where $h\dvtx {\mathcal S}_0 \to{\mathbb R}$ is defined
in~\textup{(\ref{hfun})}.

For any $\hat x=(x,\gamma) \in{\mathfrak S},$
\begin{eqnarray*}\int_{\mathfrak S}\widehat H_\kappa
(\hat x,d \hat y) \hat h( \hat y)=E_x^{-} (|\rho_0|^\kappa h(x_0)
)=\int_{{\mathcal S}_0} H_\kappa(x,d y) h(y)=\hat h(\hat x).
\end{eqnarray*}
Consequently, setting $\hat\pi_h(A \times\eta)=1/2 \pi_h(A)$ for $A
\in{\mathcal T}_0$ and $\eta\in\{-1,1\},$ we have:
\begin{eqnarray*}
&& \int_{\mathfrak S} \biggl( \int_{A \times\eta} \frac{ 1}{ \hat
h(\hat x) }\widehat H_\kappa (\hat x, d \hat y ) \hat h( \hat y)
\biggr) \hat\pi_h(d \hat x)
\\
&&\qquad =\int_{{\mathcal S}_0} \frac{ 1}{ 2 h(x) } E_x^{-} \bigl(|\rho
_0|^\kappa h(x_0); x_0 \in A \bigr) \pi_h(d x)
\\
&&\qquad =\frac{ 1}{ 2 }\int_{{\mathcal S}_0} H_\kappa(x,A) \pi_h(d
x)=\frac{ 1}{ 2 }\pi_h(A)=\hat\pi_h(A \times \eta).
\end{eqnarray*}
We will use these facts to write down formulas similar to
\textup{(\ref{expl})} for the limits $K_1(x)$ and $K_{-1}(x)$ in
\textup{(\ref{gvul1})}. Claims (b)~and~(c) of Theorem~\ref{main-rec-3}
are immediate consequences of these formulas.

First, assume that ${\mathfrak S}_1=\varnothing.$ That is, by part
(iii) of Proposition~\ref{still}, Condition G is satisfied. We get from
Proposition~\ref{gop} and \textup{(\ref{expo})} that for $\pi $-almost
every $z \in{\mathcal S}_0$ and $\eta\in\{-1,1\}$ (compare with (4.4)
in \cite{goldie}):
\begin{eqnarray*}
K_\eta(z)&=&\frac{ 1}{ 2\tilde a } \biggl[ \int_{{\mathcal
S}_0}\int_{\mathbb R}g_1(x,t)  \pi_h(dx)\,dt+\int _{{\mathcal
S}_0}\int_{\mathbb R} g_{-1}(x,t)  \pi_h(dx)\,dt \biggr]
\\
&=& \frac{ 1}{ 2\tilde a \kappa} \int_{{\mathcal S}_0} \frac {1}{h(x)}
E_x^{-} ( |R|^\kappa-|R-\xi_0|^\kappa) \pi_h(dx),
\end{eqnarray*}
where $\tilde a=\widetilde E(\log|\rho_0|).$

Assume now that ${\mathfrak S}_1 \neq\varnothing,$ that is,
Condition G is not satisfied. We get from Proposition~\ref{gop} and
\textup{(\ref {expo})} that $\pi$-a.s., if $(z,1) \in{\mathfrak
S}_\gamma$ (i.e. $z \in A_\gamma$), then
\begin{eqnarray*}
K_\eta(z)&=&\frac{ 1}{ 2\tilde a } \biggl[ \int_{A_1} \int_{\mathbb
R}g_{\eta \gamma}(x,t)  \pi_h(dx)\,dt+\int_{A_{-1}} \int_{\mathbb
R}g_{-\eta \gamma}(x,t)
 \pi_h(dx)\,dt \biggr].
\end{eqnarray*}
This completes the proof of Theorem \ref{main-rec-3}.

%s5.3 ###
\subsection[Proof of part \textup{(a)} of Theorem 1.8]{Proof of part \textup{(a)} of Theorem \textup{\protect\ref{mishne}}}
The ``if'' part of the claim is trivial. Indeed, if
\textup{(\ref{aper})} holds for a measurable function $\Gamma\dvtx
{\mathcal S}_0 \to{\mathbb R},$ then substituting
$\xi_n=\Gamma(x_{n-1})-\rho_n\Gamma(x_n)$ into the formula for $R_n$ in
\textup{(\ref{kohav11})} yields
\[
R_n= \Gamma(x_{-1}) -\Gamma(x_{n-1}) \prod_{i=0}^{n-1} \rho_i.
\]
The Markov chain induced by $(x_n)_{n \in{\mathbb Z}}$ on $({\mathcal
S}_0,{\mathcal T}_0)$ is Harris recurrent by Lemma~\ref{urg5} and hence
$P_x^{-} (|\Gamma(x_{n-1})|<M \mbox{ i.o.} ) =1$ for some $M>0.$ Since,
$P\mbox{-a.s.},$ $R_n$ converges to $R$ and $\prod_{i=0}^{n-1} \rho_i$
converges to zero, we obtain that with probability one
$R=\Gamma(x_{-1}).$ Hence for $\pi$-almost every $x \in{\mathcal S},$
$P_x^{-}(|R|>t)=0$ for all $t$ large enough.

Assume now that $\lim_{t \to\infty} t^\kappa P_x^{-}(|R|>t)=0$ for
$\pi$-almost every $x \in{\mathcal S}.$ Our aim is to show that
\textup{(\ref{aper})} holds for some measurable function $\Gamma\dvtx
{\mathcal S}\to {\mathbb R}.$ First, we will prove the following
extension of Grincevi$\check{\mbox{c}}$ius' symmetrization inequality
(cf. \cite{trakai80}, see also \cite{goldie}, Proposition~4.2 and
\cite{saporta}, Lemma 4). It will be shown in the sequel that if the
right-hand side of \textup{(\ref{sym-in})} is a.s. zero, then \textup
{(\ref{aper})} holds with the measurable function $\Gamma(x)$ defined
in~\textup{(\ref{gamd})}.

\begin{lemma}
\label{grin-sym} Let $y_n=(x_n,\xi_n,\rho_n)_{n \in{\mathbb Z}}$ be a
MMP associated with Markov chains $(x_n)_{n \in{\mathbb Z}},$ $(\xi
_n,\rho_n) \in {\mathbb R}^2,$ and let $R$ be the random variable
defined in \textup {(\ref{ar3})}. Further, for any $x \in{\mathcal S},$
let
%
%e5.3 ###
\begin{equation}\label{gamd}
\Gamma(x)=\inf \bigl\{a \in{\mathbb R}\dvtx  P_x^{-}(R \leq
a)>\tfrac{1}{2} \bigr\}.
\end{equation}
Then, for any $t>0$ and $z \in{\mathcal S},$
%
%e5.4 ###
\begin{equation}\label{sym-in}
P_z^{-}(|R| \geq t) \geq \tfrac{ 1}{ 2 } P_z^{-}\bigl(|R_n
+\Gamma(x_{n-1}) \Pi_n|>t \mbox{ for some } n \geq0\bigr),
\end{equation}
where the random variables $\Pi_n$ and $R_n$ are defined in
\textup{(\ref{kohav})} and \textup{(\ref{kohav11})}, respectively.
\end{lemma}

\begin{pf}
By its definition, $\Gamma(x)$ is a median of the random variable $R$
under the measure $P_x^{-},$ that is $P_x^{-} (R \geq \Gamma(x)
)\geq1/2$ and $P_x^{-} (R \leq\Gamma(x) ) \geq1/2.$ Moreover,
$\Gamma(x)$ is a measurable function of $x.$

Fix now any $t>0$ and let $\tau_1=\inf\{n>0\dvtx  R_n +\Gamma(x_{n-1})
\Pi_n>t\}.$ Since $\Gamma(x)$ is a median of the distribution
$P_x^{-}(R \in\cdot),$ it follows from the definition \textup{(\ref
{kohav11})} of the random variables $R_n$ and the Markov property that
\begin{eqnarray*}
P_z^{-}(R \geq t) &\geq&  \sum_{n=0}^\infty \int_{\mathcal S}P_z^{-}(
\tau_1=n; x_{n-1} \in dx,\Pi_n>0)P_x^{-}\bigl(R \geq\Gamma(x)\bigr)
\\
&&{} + \sum_{n=0}^\infty\int_{\mathcal S}P_z^{-}( \tau_1=n; x_{n-1} \in
dx,\Pi_n<0)P_x^{-}\bigl(R \leq\Gamma(x)\bigr)
\\
& \geq& \tfrac{ 1}{ 2 } P_z^{-}(\tau_1<\infty).
\end{eqnarray*}
Replacing the sequence $\xi_n$ by the sequence $-\xi_n$ and
consequently $R$ by $-R,$ we obtain [note that we can replace
$\Gamma(x_n)$ by $-\Gamma(x_n)$ because the latter is a median of
$-R$]:
\begin{eqnarray*}P_z^{-}(-R \geq t)
\geq\tfrac{ 1}{ 2 } P_z^{-}(\tau_2<\infty),
\end{eqnarray*}
where $\tau_2:=\inf\{n>0\dvtx  -R_n -\Gamma(x_{n-1}) \Pi_n>t\}.$
Combining together these two inequalities, we
get~\textup{(\ref{sym-in})}.
\end{pf}

We will apply this lemma to the Markov chain
$y_n^*=(x_n^*,Q_n^*,M_n^*)_{n \in{\mathbb Z}},$ defined below by a
``geometric sampling,'' rather than to $y_n=(x_n,\xi_n, \rho_n)_{n
\in{\mathbb Z}}.$ The stationary sequence $(x_n^*)_{n \geq-1}$ [it is
expanded then into the double-sided sequence $(x_n^*)_{n \in {\mathbb
Z}}$] is a random subsequence of $(x_n)_{n \geq-1}$ that forms a Markov
chain which inherits the properties of $(x_n)_{n \in{\mathbb Z}}$ and
in addition is \textit{ strongly aperiodic}, that is, Lemma~\ref{urg3}
holds for this chain with $d=m=1.$

Let $(\eta_n)_{n \geq0}$ be a sequence of i.i.d. variables independent
of $(x_n,\xi_n,\rho_n)_{n \in{\mathbb Z}}$ (defined in a probability
space enlarged if needed) such that $P(\eta_0=1)=1/2$ and
$P(\eta_0=0)=1/2,$ and define $\varrho_{-1}=-1,$
$\varrho_n=\inf\{i>\varrho_{n-1}\dvtx  \eta_i=1\},$ $n \geq0.$ Further,
for $n \geq-1$ let,
%
%e5.5 ###
\begin{eqnarray}\label{lkoh}
\nonumber x^*_n&=&x_{\varrho_n},
\\
Q_{n+1}^*&=&\xi_{{\varrho_n+1}}+\xi_{{\varrho_n+2}}
\rho_{{\varrho_n+1}}+\cdots + \xi_{{\varrho_{n+1}}}
\rho_{{\varrho_n+1}}\rho_{{\varrho_n+2}} \cdots
\varrho_{{\varrho_{n+1}+1}},
\\
\nonumber M_{n+1}^*&=& \rho_{{\varrho_n+1}}\rho_{{\varrho_n+2}} \cdots
\rho_{{\varrho_{n+1}}}.
\end{eqnarray}
The transition kernel of the Markov chain $(x^*_n)_{n \geq-1}$ is given
by
%
%e5.6 ###
\begin{equation}\label{hko}
H^*(x,\cdot)=\sum_{n=1}^\infty\bigl(\tfrac{1}{2} \bigr)^n H^n(x,\cdot).
\end{equation}
Hence, $(x^*_n)_{n \geq-1}$ is Harris recurrent on ${\mathcal S}_0$ and
its stationary distribution is $\pi.$ Moreover, the sequence
$(y_n^*)_{n \geq0}=(x_n^*,Q_n^*,M_n^*)_{n \geq0}$ is a stationary
Markov chain whose transitions depend only on the position of $x_n^*$
and
%
%e5.7 ###
\begin{equation}
\label{rqm} R=Q_0^*+\sum_{n=1}^\infty Q_n^* \prod_{i=0}^{n-1} M_i^*.
\end{equation}
Expand $(y_n^*)_{n \geq0}$ into a double-sided stationary sequence
$(y_n^*)_{n \in{\mathbb Z}}.$

The following corollary to Lemma~\ref{grin-sym} is immediate in view of
\textup{(\ref{rqm})}.

\begin{corollary}
Let Assumption~\textup{\ref{measure1}} hold. Then, for any $t>0$ and $z
\in {\mathcal S},$
%
%e5.8 ###
\begin{equation}\label{sym1}
P_z^{-}(|R| \geq t) \geq\tfrac{ 1}{ 2 } P_z^{-}\bigl(|R_n^*
+\Gamma(x_{n-1}^*) \Pi_n^*|>t \mbox{ for some } n \geq0\bigr),
\end{equation}
where $\Pi_n^*:=\prod_{i=0}^{n-1} M_i^*$ and $ R_n^*:=\sum_{i=0}^{n-1}
Q_i^* \Pi_i^*.$
\end{corollary}

Our aim now is to show that the right-hand side of \textup{(\ref
{sym1})} is bounded away from zero for $\pi$-almost every $z
\in{\mathcal S}.$ The main advantage of using the ``geometrically
sampled'' MMP $(x_n^*,Q_n^*,M_n^*)_{n \in{\mathbb Z}}$ is that studying
its one-step transitions one can obtain some information concerning all
possible transitions of the original MMP $(x_n,\xi_n, \rho_n)_{n
\in{\mathbb Z}}.$ We will use this when passing from \textup{(\ref
{fin})} to \textup{(\ref{fin1})} below.

At some stage of the proof, we shall apply Corollary~\ref{renewal1} to
the Markov chain $(x_n^*)_{n \in{\mathbb Z}}$ and the random walk
$V_n^*=\sum_{i=0}^{n-1} \log|M_i^*|$ considered under the measure
$\widetilde P$ introduced in Section~\ref{similar}. Let $h_\beta\dvtx
{\mathcal S}_0 \to(0,\infty)$ be the eigenfunction of the operator
$H_\beta$ in the space $B_b$ corresponding to the kernel defined in
\textup{(\ref {kbeta})}. This eigenfunction exists and is bounded away
from zero by Proposition~\ref{urg}, and it corresponds to the
eigenvalue $r_\beta$ which coincides with the spectral radius
$r_{{H_\beta}}$ of the operator. Let
%
%e5.9 ###
\begin{equation}\label{kbeta1} H_\beta^*(x,dy)=\sum_{n=1}^\infty
\bigl(\tfrac{1}{2} \bigr)^n H^n_\beta(x,\cdot).
\end{equation}
Then, similarly to \textup{(\ref{hbeta})}, $E_x^{-} (\prod_{i=1}^n
|M_i^*|^\beta)= H_\beta^{*n} \mathbf{1}(x)$ for any $\beta\geq0$ and $x
\in{\mathcal S}_0.$

Transition kernel $\widetilde H^*$ of $x_n^*$ under $\widetilde P$ is
given by
%
%e5.10 ###
\begin{equation}\label{impo} \widetilde H^*(x,dy)=\sum_{n=1}^\infty
\biggl(\frac{ 1}{ 2 } \biggr)^n \widetilde H^n= \frac{ 1}{ h(x)
}H_\kappa^*(x,dy)h(y),
\end{equation}
where as before $h(x)=h_\kappa(x).$ It follows from
\textup{(\ref{kbeta1})} that, as long as $r_\beta<2,$
\begin{eqnarray*}H^*_\beta h_\beta(x)=\frac{ r_\beta}{ 2-r_\beta}
h_\beta(x),
\end{eqnarray*}
and thus, as in Proposition~\ref{urg}, $r_\beta^*:=
\frac{ r_\beta}{ 2-r_\beta}$ is the spectral radius of the operator
$H_\beta^*$ in $B_b.$ In particular, $r^*_\kappa=1.$ Note
also that the invariant distribution of $\widetilde H^*$ coincides with
the invariant distribution $\pi_h$ of $\widetilde H.$

To enable in the use of Corollary~\ref{renewal1} we need the following
two lemmas which ensure that its conditions are satisfied. First, the
same proof as that of Lemma~\ref{alemma3} yields:

\begin{lemma}
\label{newex} Let Assumption~\textup{\ref{measure1}} hold. Then,
$\widetilde E( \log|M^*_1|)>0.$
\end{lemma}

In addition, we have:

\begin{lemma}\label{newex1}
Let Assumption~\textup{\ref{measure1}} hold. Then, the process
$\log|M^*_n|$ is nonarithmetic relative to the Markov chain $(x_n^*)_{n
\in{\mathbb Z}}$ with transition kernel $\widetilde H^*$ defined
in~\textup{(\ref{impo})} (in the sense of
Definition~\textup{\ref{arc}}).
\end{lemma}

\begin{pf}
Since the process $\log|\rho_i|$ is nonarithmetic relative to the
Markov chain $(x_n)_{n \in{\mathbb Z}}$ with kernel $\widetilde H,$ the
claim follows from Lemma~A.6 in~\cite{alsmeyer}, which deals with the
nonarithmetic condition relative to the ``sampled'' Markov chain
$(x_n^*)_{n \in{\mathbb Z}}.$%\rightqed
\end{pf}

We are now in position to complete the proof of part (a) of
Theorem~\ref{mishne}.

\begin{lemma}
Let Assumption~\textup{\ref{measure1}} hold and suppose in addition that\break
$\lim_{t \to \infty} t^\kappa P_x^{-}(|R|>t)=0$ for $\pi$-almost every
$x \in {\mathcal S}.$ Then, \textup{(\ref{aper})} holds with the
function $\Gamma(x)$ defined in~\textup{(\ref{gamd})}.
\end{lemma}

\begin{pf}
For $n \in{\mathbb Z},$ let $\alpha_n=R_n^* +\Gamma(x_{n-1}^*) \Pi
_n^*$ and write $\alpha_n=\alpha_{n-1}+\beta_n,$ where
\begin{eqnarray*}\beta_n &=& Q^*_{n-1}
\Pi^*_{n-1}+\Gamma(x^*_{n-1}) \Pi_n^*-\Gamma(x^*_{n-2}) \Pi^*_{n-1}
\\
&=& \Pi^*_{n-1}
\bigl(Q^*_{n-1}+\Gamma(x^*_{n-1})M^*_{n-1}-\Gamma(x^*_{n-2}) \bigr).
\end{eqnarray*}
Set
%
%e5.11 ###
\begin{equation}\label{deltan}
\delta_n:=Q^*_n+\Gamma(x^*_n)M^*_n-\Gamma(x^*_{n-1}).
\end{equation}
Thus, $\alpha_n=\alpha_{n-1}+\beta_n=\alpha_{n-1}+\Pi^*_{n-1}
\delta_{n-1},$ and hence for any $\varepsilon>0$ (cf. \cite{goldie},
page~157):
\begin{eqnarray*}
P_z^{-}(|\alpha_n| > t \mbox{ for some }n \geq 0) &\geq&
P_z^{-}(|\beta_n| > 2t \mbox{ for some }n \geq0)
\\
&\geq& P_z^{-} (|\Pi_n^*|>2t/\varepsilon \mbox{ and }
|\delta_n|>\varepsilon \mbox{ for some } n \geq1).
\end{eqnarray*}
Indeed, $|\beta_n| >2t$ implies that either $|\alpha_{n-1}|>t$ or, if
not, $|\alpha_n| \geq|\beta_n| -|\alpha_{n-1}|>2t-t=t.$

Fix a number $\varepsilon>0$ and let $\upsilon(t)=\inf\{n \geq1\dvtx
|\Pi_n^*|>2t/\varepsilon\}.$ Then, setting
%
%e5.12 ###
\begin{equation}
\label{mrw} V_n^*:=\log|\Pi_n^*|=\sum_{i=0}^{n-1} \log|M_i^*|,
\end{equation}
we obtain from \textup{(\ref{sym1})} and the Markov property that for
any $z \in {\mathcal S}_0,$
\begin{eqnarray*}
&& t^\kappa P_z^{-}(|R| \geq t)
\\
&&\qquad \geq \frac{ t^\kappa }{ 2 } \int_{{\mathcal S}_0}
P_z^{-}\bigl(x^*_{\upsilon(t)-1} \in dx,
|\delta_{\upsilon(t)}|>\varepsilon,\upsilon(t)<\infty\bigr)
\\
&&\qquad =  \frac{ t^\kappa}{ 2 } E_z^{-}
\bigl(P^{-}_{x^*_{\upsilon(t)-1}}(|\delta_0|>\varepsilon);
\upsilon(t)<\infty\bigr)
\\
&&\qquad = \frac{ 1}{ 2 } \biggl(\frac{ \varepsilon}{ 2 } \biggr)^\kappa h(z)
\widetilde E^{-}_{z}
\bigl(e^{-\kappa(V_{\upsilon(t)}^*-\log(2t/\varepsilon))}
P^{-}_{x_{\upsilon(t)-1}^*}(|\delta_0|>\varepsilon)/
h\bigl(x_{\upsilon(t)-1}^*\bigr) \bigr),
\end{eqnarray*}
where the expectation $\widetilde E_z^{-}$ is according to the measure
$\widetilde P_z^{-}$ defined in Section~\ref{similar}.

Thus, in virtue of part (b) of Theorem~\ref{main-rec-3} it suffices to
prove that under Assumption~\ref{measure1},
\begin{longlist}[(ii)]
\item[] either

\item[(i)] for some $\varepsilon>0$ and probability measure $\hat\pi$
absolutely continuous with respect to $\pi,$ either the following limit
exists and is strictly positive:
%
%e5.13 ###
\begin{equation}\label{str}
\lim_{t \to\infty} \widetilde E_{\hat\pi}^- \bigl(e^{-\kappa
(V_{\upsilon(t)}^*-\log(2t/\varepsilon))} P^{-}_{x_{\upsilon(t)-1}^*}
(|\delta_0|>\varepsilon) \bigr),
\end{equation}
where $\widetilde E^-_{\hat\pi}(\cdot):= \int_{{\mathcal S}_0}
\widetilde E_z^-(\cdot)\hat\pi(d z),$

\item[] or, if not,

\item[(ii)] then, \textup{(\ref{aper})} holds with the function
$\Gamma(x)$ defined in \textup{(\ref{gamd})}.
\end{longlist}

To bound the limit in \textup{(\ref{str})} away from zero we will apply
Corollary~\ref{renewal1} to the Markov chain $(x^*_n)_{n \in{\mathbb
Z}}$ on $({\mathcal S}_0,{\mathcal T}_0)$ introduced in
\textup{(\ref{lkoh})} and governed by the kernel $\widetilde H^*$
defined in \textup{(\ref{impo})}, the random walk $V_n^*$ defined in
\textup{(\ref{mrw})}, and the function
\[
g(x,t)=e^{- \kappa t} P_x^{-} (|\delta_0|>\varepsilon ).
\]

Let $\sigma_{-1}=-1,$ $V^*_{-1}=0,$ and for $n \geq0,$
$\sigma_n=\inf\{i>\sigma_{n-1}\dvtx  V^*_i> V^*_{\sigma_{n-1}}\}.$
Further, let $\hat\pi$ be the stationary distribution of the Markov
chain $\hat x_n:=x_{\sigma_n}^*$ under $\widetilde P$ (which exists and
is unique by \cite{alsmeyer}, Theorem~4). The measure $\hat \pi$ is an
irreducible measure of the Markov chain $(x_n^*)_{n \in {\mathbb Z}}$
with transition kernel $\widetilde H^*$ and hence is absolutely
continuous with respect to its stationary distribution, which in turn
is equivalent to the stationary distribution $\pi$ of $(x_n^*)_{n \in
{\mathbb Z}}$ with transition kernel $H^*.$

To apply Corollary~\ref{renewal1} to the Markov chain $(x^*_n)_{n \in
{\mathbb Z}}$ with kernel $\widetilde H^*$ and the random walk $V^*_n,$
we need to check conditions \textup{(\ref{al1})} and
\textup{(\ref{al5})} for the function
\[
b(x,y)=\widetilde E_x^- \bigl( e^{- \kappa(\widehat V_0 -y)}
\mathbf{1}_{\{ \widehat V_0>y\}} P_{\hat x_0}^{-}
(|\delta_0|>\varepsilon) \bigr),
\]
where $\widehat V_n:=V_{\sigma_n}^*.$

Condition \textup{(\ref{al1})} follows from the following estimate
valid for any $\delta>0$:
\begin{eqnarray*}
|b (x,y+ \delta)-b(x,y)| &\leq& \widetilde E^-_x \bigl(
\big|e^{\kappa\delta} \mathbf{1}_{\{\widehat V_0>y+\delta
\}}-\mathbf{1}_{\{\widehat V_0>y\}}\big| \bigr)
\\
&=&  (e^{\kappa\delta}-1) \widetilde P_x^- (\widehat V_0>y +\delta)+
\widetilde P_x^- ( y \leq\widehat V_0 < y+\delta).
\end{eqnarray*}
As to condition \textup{(\ref{al5})}, we have:
\[
b(x,y) \leq \cases{ e^{\kappa y}, &\quad if $y<0$, \cr \widetilde E^-_x
\bigl( \mathbf{1}_{\{\widehat V_0>y\}} \bigr)= \widetilde P^-_x (\widehat V_0>y
), &\quad if $y \geq0$.}
\]
Hence,
\begin{eqnarray*}
&& \int_{{\mathcal S}_0} \sum_{n \in{\mathbb Z}} \sup_{n \leq y< n+1 }
|b(x,y)| \hat\pi(dy)
\\
&&\qquad \leq \sum_{n=0}^\infty e^{-\kappa n} + \int_{{\mathcal
S}_0}\sum_{n=0}^\infty\widetilde P_x^- (\widehat V_0>n )\hat \pi(dx)<
\infty,
\end{eqnarray*}
because by part (iv) of \cite{alsmeyer}, Theorem~2,
\[
\int_{{\mathcal S}_0} \sum_{n=0}^\infty\widetilde P^-_x (\widehat V_0>n
)\hat\pi(dx) \leq\int_{{\mathcal S}_0} \widetilde E_x^- (\widehat V_0
)\hat\pi(dx)<\infty.
\]
(Part (iv) of \cite{alsmeyer}, Theorem~2 implies that the constant
$\mu_1$ in the statement of Corollary~\ref{renewal1} is finite. In our
case, $\mu_1= \int _{{\mathcal S}_0} \widetilde E^-_x (\widehat V_0
)\hat\pi(dx)$.)

Let $\widehat H$ be the transition kernel of the Markov chain $(\hat
x_n, \widehat V_n-\widehat V_{n-1})_{n \geq0}.$ It follows from
Corollary~\ref{renewal1} that for some $A \in(0,1),$
\begin{eqnarray*}
&&\lim_{t \to\infty} \widetilde E_{\hat\pi}^- \bigl(e^{-\kappa
(V_{\upsilon(t)}^*-\log(2t/\varepsilon))}
P^{-}_{x_{\upsilon(t)}^*-1}(|\delta_0|>\varepsilon) \bigr)
\\
&&\qquad = \frac{1}{\mu_1} \int_{{\mathcal S}_0} \int_{{\mathcal S}_0
\times (0,\infty)} \int_{[0,z)} e^{-\kappa w}
P^{-}_y(\delta_0>\varepsilon)\,d w\, \widehat H(x,dy \times dz)
\hat\pi(dx)
\\
&&\qquad  = \frac{1}{\mu_1 \kappa} \int_{{\mathcal S}_0}
\int_{{\mathcal S}_0 \times (0,\infty)} (1-e^{-\kappa z})
P^{-}_y(|\delta_0|>\varepsilon) \widehat H(x,dy \times dz) \hat\pi(dx)
\\
&&\qquad \geq A \int_{{\mathcal S}_0} \int_{{\mathcal S}_0}
P^{-}_y(|\delta _0|>\varepsilon) \widehat H \bigl(x,dy \times(0,\infty)
\bigr) \hat\pi(dx)
\\
&&\qquad = A\int_{{\mathcal S}_0} P^{-}_y(|\delta_0|>\varepsilon)
\hat\pi(dy) =A P^{-}_{\hat\pi}(|\delta_0|>\varepsilon).
\end{eqnarray*}
It follows that if \textup{(\ref{str})} is not true for any
$\varepsilon>0$ then
%
%e5.14 ###
\begin{equation}
\label{fin} P^{-}_{\hat\pi}(\delta_0=0)=1.
\end{equation}
It remains to show that \textup{(\ref{fin})} implies that
\textup{(\ref{aper})} holds for the function $\Gamma$ defined
in~\textup{(\ref{gamd})}. By the definition of the kernel $H^*$ in
\textup{(\ref{hko})} and the quantity $\delta_n$ in
\textup{(\ref{deltan})}, we get from \textup{(\ref {fin})} that
%
%e5.15 ###
\begin{equation}
\label{fin1} P^{-}_{\hat\pi} \bigl(R_n+\Gamma(x_n) \Pi_n-\Gamma
(x_{-1}) =0 \bigr)=1\qquad\mbox{for all } n \in{\mathbb N}.
\end{equation}
Taking respectively $n=0$ and $n=1$ in the last equality we obtain that
$P^{-}_{\hat\pi} (\xi_0+\Gamma(x_0) \rho_0-\Gamma(x_{-1}) =0 )=
P^{-}_{\hat\pi} (\xi_0+\xi_1\rho_0+\Gamma(x_1)
\rho_0\rho_1-\Gamma(x_{-1}) =0 )=1.$ It follows that
\[
P^{-}_{\hat\pi} \bigl(\xi_1+\Gamma(x_1) \rho_1-\Gamma(x_0) =0 \bigr)=1.
\]
Similarly, by induction on $n,$ one can show that
\[
P^{-}_{\hat\pi} \bigl(\xi_n+\Gamma(x_n) \rho_n-\Gamma(x_{n-1}) =0
\bigr)=1\qquad \mbox{for all } n \in{\mathbb N}.
\]
Since the Markov chain $(x_n)$ is $\pi$-recurrent and $\hat\pi$ is
absolutely continuous with respect to $\pi,$ we
obtain~\textup{(\ref{aper})}.
\end{pf}

%s5.4 ###
\subsection[Proof of parts \textup{(b)} and \textup{(c)} of Theorem 1.8]{Proof
of parts \textup{(b)} and \textup{(c)} of Theorem \textup{\protect\ref{mishne}}}
Let ${\mathcal S}_0$ be as defined in Lemma~\ref{urg3} and recall the
regeneration times $N_n$ defined in Section~\ref{bagr}. Let
$Q_0=\xi_0+\mathbf{1}_{\{N_1 \geq1\}} \sum_{i=0}^{N_1-1}\xi_{i+1}
\prod_{j=0}^i \rho_j$ and $M_0=\prod_{i=0}^{N_1} \rho_i,$ and for $n
\geq1,$
\begin{eqnarray*}
Q_n&=& \xi_{N_n+1}+\mathbf{1}_{\{N_{n+1} -N_n \geq2\}}
\sum_{i=N_n+1}^{N_{n+1}-1} \xi_{i+1}\prod_{j=N_n+1}^i
\rho_j\quad\mbox{and}\quad M_n=\prod_{i=N_n+1}^{N_{n+1}} \rho_i.
\end{eqnarray*}
The pairs $(Q_n,M_n),n \geq0,$ are one-dependent and for $n \geq1$ they
are identically distributed. Since the series in \textup{(\ref{ar3})}
converges absolutely, we obtain the representation
%
%e5.16 ###
\begin{equation}\label{arko}
R=Q_0+M_0 \bigl(Q_1+M_1 \bigl(Q_2+M_2 (Q_3+\cdots) \bigr)
\bigr):=Q_0+M_0 \widehat R.
\end{equation}
Note that $x_{{N_1}}$ is distributed according to the measure $\psi$
introduced in Lemma~\ref{urg} and hence $P(|\widehat
R|>t)=P_\psi^{-}(|R|>t),$ where we denote as usual
$P^{-}_\psi(\cdot):=\int_{\mathcal S}P^{-}_x(\cdot)\psi(d x).$ We have:

\begin{lemma}
\label{psi-gvul} The following limit exists and is strictly positive:
%
%e5.17 ###
\begin{equation}\label{psig}
\widehat K=\lim_{t \to\infty} t^\kappa P^{-}_\psi(|R|>t)=\lim_{t
\to\infty} t^\kappa P(|\widehat R|>t).
\end{equation}
\end{lemma}

\begin{pf}
The measure $\psi$ is an irreducible measure of the Markov chain
$(x_n)_{n \in{\mathbb Z}}$ and hence it is absolutely continuous with
respect to its stationary distribution $\pi.$ Therefore, the claim
follows by the bounded convergence theorem from part~(a) of
Theorem~\ref{main-rec-3} and part~(c) of Lemma~\ref{alemma1}.
\end{pf}

We will show next that the contribution of $Q_0$ in $R$ is negligible
in the following precise sense [recall that $\xi_n$ are assumed to be
bounded by~(\ref{calso})]: for some $\beta>\kappa,$
%
%e5.18 ###
\begin{equation}\label{233}
\sup_{x \in{\mathcal S}_0} E_x^{-} \Biggl( \Biggl[\mathbf{1}_{\{N_1
\geq 1\}}\sum_{i=0}^{N_1-1} \prod_{j=0}^i |\rho_j|
\Biggr]^\beta\Biggr)<\infty.
\end{equation}
Let $A(x)=E_x^{-} ( [\mathbf{1}_{\{N_1 \geq1\}}\sum_{i=0}^{N_1-1}
\prod_{j=0}^i |\rho_j| ]^\beta)<\infty.$ Since for any positive numbers
$\{a_i\}_{i=1}^n$ we have $(a_1+a_2+\cdots + a_n)^\beta\leq n^\beta
(a_1^\beta+a_2^\beta+\cdots +a_n^\beta),$ we obtain for any $\beta>0$
and $x \in{\mathcal S}_0$:
%
%e5.19 ###
\begin{eqnarray}\label{s33}
\nonumber A(x)&=&  E_x^{-} \Biggl(\,\sum_{n=1}^\infty\sum_{i=0}^{n-1}
\prod_{j=0}^i |\rho_j| \mathbf{1}_{\{N_1=n\}} \Biggr)^\beta
\\
&=&  \sum_{n=1}^\infty E_x^{-} \Biggl(\, \sum_{i=0}^{n-1} \prod_{j=0}^i
|\rho_j| \mathbf{1}_{\{N_1=n\}} \Biggr)^\beta
\\
\nonumber &\leq& \sum_{n=1}^\infty n^\beta\sum_{i=0}^{n-1} E_x^{-}
\Biggl(\,\prod_{j=0}^i |\rho_j|^\beta\mathbf{1}_{\{N_1 \geq n\}}
\Biggr).
\end{eqnarray}
Let
\[
\widetilde \Theta_\beta(x,dy):= \Theta(x,dy) E (|\rho_0\rho_1 \rho_2
\cdots \rho_{m-1}|^\beta| x_{-1}=x,x_{m-1}=y ),
\]
where the kernel $\Theta(x,d y)$ on $({\mathcal S}_0,{\mathcal T}_0)$
is defined in \textup{(\ref{theta})}, and let
\begin{eqnarray*}
K_\beta(x,dy)&:=&H^m(x,dy)E (|\rho_0\rho_1 \rho_2 \cdots
\rho_{m-1}|^\beta|x_{-1}=x,x_{m-1}=y )
\\
&\hspace*{3pt} = &H^m_\beta(x,dy),
\end{eqnarray*}
where the kernel $H_\beta$ on $({\mathcal S}_0,{\mathcal T}_0)$ is
defined in \textup{(\ref{kbeta})}.

Then for any $x \in{\mathcal S}_0,$
\[
\widetilde\Theta_\beta\mathbf{1}(x)=E_x^{-} \Biggl(\prod_{j=0}^{m-1}
|\rho_j|^\beta\mathbf{1}_{\{N_1 \geq m\}} \Biggr)\quad\mbox{and}\quad
K_\beta\mathbf{1}(x)=E_x^{-} \Biggl(\prod_{j=0}^{m-1}
|\rho_j|^\beta\Biggr).
\]
By Lemma~\ref{urg3} and \textup{(\ref{calso})}, the kernels $K_\beta$
and $\widetilde\Theta_\beta$ satisfy the conditions of Proposition
\ref{urg} with $ s(x,y)=E (|\rho_0\rho_1 \rho_2 \cdots
\rho_{m-1}|^\beta| x_{-1}=x,x_{m-1}=y )$ and $c_1=c_\rho^{-m}.$ In
virtue of Lemma~\ref{urg5}, the spectral radius of $H_\kappa$ and hence
$K_\kappa$ is equal to 1. Thus, by part~(c) of Proposition \ref{urg},
the spectral radius of $\widetilde\Theta_\kappa$ is strictly less than
one. Since $r_{\widetilde\Theta_\beta}$ is a continuous function of
$\beta,$ we have for some $\beta>\kappa$:
%
%e5.20 ###
\begin{equation}\label{betak1}
r_{\widetilde\Theta_\beta}<1.
\end{equation}
For $l \in{\mathbb N},$ denote $ \hat l=m \cdot\max\{[l/m],1\},$ where
$m$ is as in \textup{(\ref{lem-lem})}. We obtain from
\textup{(\ref{betak1})} that for any $l \in{\mathbb N}, n>\max
\{l,m\},$ $x \in {\mathcal S}_1,$ and for suitable constants
$A_\beta>0$, $\widetilde \Lambda_\beta<0$:
\begin{eqnarray*}
E_x^{-} \Biggl(\,\prod_{j=0}^l |\rho_j|^\beta\mathbf{1}_{\{ N_1 \geq
n\}} \Biggr) &\leq& c_\rho^m E_x^{-} \Biggl(\,\prod_{j=0}^{\hat l-1}
|\rho_j|^\beta\mathbf {1}_{\{N_1 \geq \hat n\}} \Biggr)
\\
&\leq& c_\rho^m \widetilde\Theta _\beta^{\hat l/m} \Theta^{(\hat n-\hat
l)/m} \mathbf{1}(x) \leq A_\beta e^{n \widetilde \Lambda_\beta},
\end{eqnarray*}
where in the first inequality we use \textup{(\ref{calso})} and the
fact that $\hat n \leq n$ for $n>m$ (note also that $r_\theta<1$ by
Proposition~\ref{urg} applied to the kernels $H$ and $\Theta$). This
yields \textup{(\ref{233})} in virtue of \textup{(\ref{s33})}.

Fix some $\beta>\kappa$ which satisfies \textup{(\ref{233})} and
$\alpha\in (\frac{\kappa}{\beta}, 1 ).$ By \textup{(\ref{s33})} and the
Chebyshev inequality, $\lim_{t \to\infty}t^\kappa P_x^{-}(|Q_0| \geq
t^\alpha)=0$ uniformly in $x.$ Let
\begin{eqnarray*}
M_{0,1}=\mathbf{1}_{\{N_1-m=-1\}}+\mathbf{1}_{\{N_1-m \geq0\}} \cdot
\prod_{i=0}^{N_1-m} |\rho_i|\quad \mbox{and}\quad  M_{0,2}=
\prod_{i=N_1-m+1}^{N_1} |\rho_i|.
\end{eqnarray*}
Then, $M_0=M_{0,1}\cdot M_{0,2}$ and $c_\rho^{-m} M_0 \leq M_{0,1} \leq
c_\rho^m M_0,$ where $c_\rho$ is introduced in assumption~(A4).

Recall the random variable $\widehat R$ defined in
\textup{(\ref{arko})} and note that $M_{0,1}$ and $\widehat R$ are
independent under the measure $P^{-}_x$ because only the $m-1$ last
variables in the block $(x_0,x_1,\ldots,x_{{N_1-1}}\}$ are dependent on
$x_{{N_1}}.$

For any $\beta>\kappa$ such that \textup{(\ref{233})} holds, we have
\begin{eqnarray*}
t^\kappa P^{-}_x( |R|>t) &\leq& t^\kappa P^{-}_x (|Q_0|+|M_0\widehat
R|>t, |Q_0|<t^\alpha)+t^\kappa P^{-}_x (|Q_0| \geq t^\alpha)
\\
&\leq& t^\kappa P^{-}_x (|M_0\widehat R|>t-t^\alpha )+\frac{ t^\kappa}{
t^{\alpha\beta} } E_x^{-} (|Q_0|^\beta)
\\
&\leq& t^\kappa P^{-}_x (c_\rho^m |M_{0,1} \widehat
R|>t-t^\alpha)+E_x^{-} (|Q_0|^\beta).
\end{eqnarray*}
The expectation $E_x^{-} (|Q_0|^\beta)$ is bounded on ${\mathcal S}_0$
by \textup{(\ref{233})}, while \textup{(\ref{psig})} and the fact that
$\widehat R$ is independent of $M_{0,1}$ under $P_x^{-}$ imply that for
some $L>0,$
\begin{eqnarray*}
t^\kappa P^{-}_x (c_\rho^m| M_{0,1} \widehat R|>t-t^\alpha) \leq L
\biggl( \frac{ t}{ t-t^\alpha} \biggr)^\kappa E_x^{-} (
|M_{0,1}|^\kappa)\qquad \forall t>1
\end{eqnarray*}
yielding the upper bound in \textup{(\ref{bound})} since the
expectation $E_x^{-} |(M_{0,1}|^\beta)$ is bounded on ${\mathcal S}_0$
in view of \textup{(\ref{233})}.

To get the lower bound in \textup{(\ref{43-bound-43})}, write
\begin{eqnarray*}
t^\kappa P^{-}_x( |R|>t) &\geq& t^\kappa P^{-}_x (|M_0\widehat
R|-|Q_0|>t )
\\
&\geq& t^\kappa P^{-}_x (|M_0\widehat R|-|Q_0|>t, |Q_0|<t^\alpha)
\\
&\geq& t^\kappa P^{-}_x (|M_0\widehat R|>t+t^\alpha)-P^{-}_x
(|Q_0|>t^\alpha)
\\
&\geq& t^\kappa P^{-}_x (|M_0\widehat R|>t+t^\alpha)-\frac{ t^\kappa}{
t^{\alpha \beta} } E_x^{-}(|Q_0|^\beta)
\\
&\geq& t^\kappa P^{-}_x (c_\rho^{-m} |M_{0,1}\widehat R|>t+t^\alpha
)-\frac{ t^\kappa}{ t^{\alpha\beta} } E_x^{-} (|Q_0|^\beta),
\end{eqnarray*}
and note that $\frac{ t^\kappa}{ t^{\alpha\beta} } E_x^{-}
(|Q_0|^\beta)$ converges to zero uniformly on $x$ by
\textup{(\ref{233})} while by \textup{(\ref{psig})} we have for any
$\lambda>0,$ some constant $J>0$ that depends on $\lambda,$ and all $t$
large enough:
\begin{eqnarray*}
t^\kappa P^{-}_x (c_\rho^{-m} |M_{0,1} \widehat R|>t ) &\geq& t^\kappa
P^{-}_x (\lambda\cdot c_\rho^{-m} \cdot|\widehat R|>t;|
M_{0,1}|>\lambda)
\\
&\geq& J P^{-}_x (|M_{0,1}| \geq\lambda).
\end{eqnarray*}
To complete the proof it remains to show that for some $\lambda>0$
there exists a number $\delta_1>0$ such that
\begin{eqnarray*}P^{-}_x (|M_{0,1}|
\geq\lambda)>\delta_1,\qquad\pi\mbox{-a.s.}
\end{eqnarray*}
Toward this end observe that for every $x \in{\mathcal S}_0,$ with
$\vartheta\in{\mathbb N}$ defined in \textup{(\ref{as})} and $c_\rho>0$
defined in~(A4) (we will assume, actually without loss of generality,
that $\vartheta>m$),
\begin{eqnarray*}
P^{-}_x \bigl(|M_{0,1}| \geq c_\rho^{-(\vartheta-m)} \bigr) &\geq&
P^{-}_x \biggl( |M_{0,1}| \geq\min_{m \leq i
\leq\vartheta}c_\rho^{-(i-m)};N_1 \leq \vartheta \biggr)
\\
&=& P^{-}_x (N_1 \leq\vartheta) \geq\delta,
\end{eqnarray*}
where $\delta>0$ is defined in \textup{(\ref{as})}.

\begin{appendix}%\appendix
%s6 ###
\section{\texorpdfstring{Proof of Proposition \protect\lowercase{\ref{urg}}}{Appendix A: Proof of Proposition 2.4}} \label{proofp}
(a) First, we
note that if a nonnegative eigenfunction $f \not\equiv0$ of the
operator $K\dvtx  B_b \to B_b$ exists then necessarily $\inf_x f(x)>0.$
Indeed, assuming that $K f=\lambda f$ for some $\lambda>0,$ we have for
any $x \in{\mathcal S}_0,$
\begin{eqnarray*}
\sum_{i=1}^{d+m} \lambda^i f(x) &=& \sum_{i=1} ^{d+m} K^i f(x) \geq
\sum_{i=1} ^d \int_{{\mathcal S}_1} \int_{{\mathcal S}_1} K^i(x,dz)
K^m(z,dy)f(y)
\\
&\geq& p \cdot c_1^{-1} \cdot\int_{{\mathcal S}_1} f(y)\psi(dy)>0 ,
\end{eqnarray*}
where the last inequality follows from the fact that $f(x)>0$ for every
$x \in{\mathcal S}_0$ (cf.~\cite{nummelin}, Proposition 5.1(ii)).

The proof of the existence of such $f \in B_b$ is an application of
Nussbaum's extension of the Krein--Rutman theorem
(cf. Theorem~2.2 in~\cite{nussbaum}).
Theorem~2.2. In view of this theorem (this is explained in
Appendix~\ref{nuss}) it is sufficient to show that there exists a
double-indexed sequences of compact linear operators $Q_{n,i}$ on the
space $B_b$ such that
%
%e6.1 ###
\begin{equation}\label{version}
\limsup_{i \to \infty} \sqrt[i]{\rule{0pt}{7pt}\smash{\|K^i-Q_{n,i} \|}} \leq1/n,\qquad n
\in{\mathbb N}.
\end{equation}
It even suffices to show that $\limsup_{i \to\infty} \sqrt[i]{\|K^{m
i}-\widehat Q_{n,i} \|} \leq1/n$ for some compact operators $\widehat
Q_{n,i}$ on $B_b,$ since we can then set $Q_{n,i}=K^{i-m j_i}\widehat
Q_{n,j_i},$ where $j_i$ is the integer part of $i/m.$ For this purpose
we shall adapt the Yosida--Kakutani's proof that Markov kernels
satisfying Doeblin's condition are quasi-compact (cf. \cite{yosida},
Section~4.7).

\begin{longlist}[(1)]
\item[(1)] First, we observe that if $n(x,y)$ and $j(x,y)$ are jointly
measurable bounded function, then the product of the two operators
defined by the kernels $N(x,d y)=n(x,y)\mu(d y)$ and $J(x,d
y)=j(x,y)\mu(d y)$ is compact in $B_b.$ Indeed, we can approximate
$n(x,y)$ in $L_1({\mathcal S}_0 \times{\mathcal S}_0,{\mathcal T}_0
\times {\mathcal T}_0,\mu\times\mu)$ up to $1/i$ by a simple function
$n_i(x,y)$ which is a finite linear combination of the indicator
functions of ``rectangle'' sets $B_{i,k} \times C_{i,k},$ where
$B_{i,k},C_{i,k} \subset{\mathcal S}_0.$ Then, the operators
corresponding to the kernels $N_i(x,d y)=n_i(x,y) \mu(d y)$ are
finite-dimensional and hence $J N=\lim_{i \to\infty} J N_i,$ being the
limit in operator norm of a sequence of compact operators, is compact.

\item[(2)] Fix $n \in{\mathbb N}$ and let $\delta=\delta(1/n)$ be
defined as in condition (iii) of the proposition. Let $k(x,y)$ be a
jointly measurable density of the kernel $K^m$ with respect to $\mu$
(such a density exists since the $\sigma$-field ${\mathcal T}_0$ is
assumed to be countably generated, see, e.g., \cite{nummelin},
Lemma~2.5) and set
\[
q_n(x,y)= \min\{k(x,y), \delta^{-1} \cdot\|K^m\|\}.
\]
Let $D_x=\{y \in{\mathcal S}_0\dvtx  k(x,y) \neq q_n(x,y)\},$ thus
$k(x,\cdot) \geq\delta^{-1} \|K^m\|$ on $D_x.$ Since
\[
\sup_x K^m(x,D_x)=\sup_x \int_{D_x} k(x,y) \mu(dy) \leq \|K^m\|,
\]
then $\mu(D_x) \leq\delta.$ Hence, letting $Q_n(x,d y)=q_n(x,d y) \mu(d
y),$
\[
\|K^m-Q_n\| \leq\sup_x \int_{D_x} k(x,y) \mu(dy)= \sup_x K^m(x,D_x)
\leq1/n.
\]

\item[(3)] Let $R_n=K^m-Q_n.$ Then $K^{m i}=(Q_n+R_n)^i=\sum$ of $2^i$
terms each of them, except maybe those $i+1$ where $Q_n$ appear at most
once, is compact by~(1). But
\begin{eqnarray*}
&& \|R_n^i+Q_n R_n^{i-1}+R_nQ_nR_n^{i-2}+\cdots+R_n^{i-1}Q_n\|
\\
&&\qquad \leq (1/n)^i+i \cdot\|Q_n\| \cdot(1/n)^{i-1} \leq(1/n)^i+i
\cdot\|K^m\| \cdot (1/n)^{i-1},
\end{eqnarray*}
as required.

(b) The proof for the kernel $\widehat\Theta$ on $({\mathcal
S}_1,{\mathcal T}_1)$ is the same as for $K,$ since the conditions of
this proposition hold for $\widehat \Theta$ as well (with $d=m=1$).

(c) Let $c_{_K}>1$ be a constant such that $f(x) \in
(c_{_K}^{-1},c_{_K})$ for all $x \in{\mathcal S}_0.$ Then, for any $x
\in {\mathcal S}_0,$ $c_{_K}^{-1} f(x) \leq\mathbf{1}(x) \leq c_{_K}
f(x),$ and hence
%
%e6.2 ###
\begin{equation}\label{bbb}
c_{_K}^{-2} r_{_K}^n \leq K^n \mathbf{1}(x) \leq c_{_K}^2
r_{_K}^n\qquad \forall x \in{\mathcal S}_0.
\end{equation}
Let $\widehat K(x,\cdot)$ be the restriction of the kernel $K^m$ to the
states of the set ${\mathcal S}_1.$ It follows from
\textup{(\ref{bbb})} that the spectral radius of $\widehat K$ coincides
with $r_{_K}^m.$

By \cite{nummelin}, Proposition~5.3 and \cite{nummelin}, Theorem~5.2,
the kernel $\widehat\Theta$ has an invariant measure
$\pi_{\widehat\Theta}.$ Since $r_{_K}^m f \geq \widehat\Theta f,$ the
equality $r_{\widehat\Theta}=r_{_K}^m$ would imply by \cite{nummelin},
Proposition~5.3 and \cite{nummelin}, Theorem~5.1 that
$\pi_{\widehat\Theta}\mbox{-a.s.},$ $\widehat \Theta f (x)=r_{_K}^m
f(x) =K^m f(x) ,$ which is impossible because $f(x)>0$ and $K^m(x,d
y)-\widehat\Theta(x,d y) \geq r c_1^{-1} \psi(d y)$ for any $x
\in{\mathcal S}_1.$ Hence $r_{\widehat\Theta}< r_{_K}^m.$
\end{longlist}

%s7 ###
\section{\texorpdfstring{The Nussbaum fixed point theorem}{Appendix B: The Nussbaum fixed point theorem}}\label{nuss}
This appendix is devoted to the Nussbaum's extension of the
Krein--Rutman fixed point theorem (cf. Theorem~2.2 in~\cite{nussbaum}) or, to
be precise, to the version of this theorem which is actually used
in~\textup{(\ref{version})}.

Let $X$ be a Banach space. For a bounded subset $S$ of $X,$
Kuratowski's measure of noncompactness $\alpha(S)$ is defined by
\begin{eqnarray*}
\alpha(S)=\inf\Biggl\{d>0\dvtx  S=\bigcup_{i=1}^n S_i, n \in{\mathbb N},
\mbox{ and } D(S_i) \leq d \mbox{ for } 1 \leq i \leq n\Biggr\},
\end{eqnarray*}
where $D(S):=\sup_{x,y \in S}\|x-y\|$ is the diameter of the set~$S$.

A bounded linear operator $K$ in $X$ is called a
\textit{$b$-set-contraction} for a number $b \geq0$ if $ \alpha(K(S))
\leq b \alpha(S)$ for every bounded subset $S$ of $X.$ A closed subset
$C$ of $X$ is called a \textit{ cone} if the following holds: (i) if
$x,y \in C$ and $\alpha,\beta\geq0$ are nonnegative reals, then $\alpha
x+\beta y \in C.$ (ii) if $x \in C-\{0\},$ then $-x \notin C.$

\begin{theorem}[\textup{(\cite{nussbaum}, Theorem~2.2)}]\label{nuss1}
Let $X$ be a Banach space, $C$ be a cone in $X,$ and $K$ be a bounded
linear operator in $X$ such that $K(C) \subset C.$ Let
\[
\|K\|_{_C}:=\sup\{\|Ku\|\dvtx  u \in C,\|u\| \leq1\}
\]
and $\alpha_{_C}(K):= \inf\{b\geq0\dvtx  K_{_C} \mbox{ is a} b\mbox
{-set-contraction}\},$ where $K_{_C}\dvtx C \to C$ is the restriction
of $K$ to the cone $C.$ Further, let
\[
r_{_C}(K):=\lim_{n \to\infty}
\sqrt[n]{\rule{0pt}{7pt}\smash{\|K^n\|_{_C}}}
\quad\mbox{and}\quad \rho_{_C}(K):=\lim_{n \to\infty}e
\sqrt[n]{\rule{0pt}{7pt}\smash{\alpha_{_C}(K^n)}}.
\]
Assume that $\rho_{_C}(K)< r_{_C}(K).$ Then there exists an $x \in C-\{
0\}$ such that $K x=r_{_C}(K) x.$
\end{theorem}

We want to apply this theorem in the situation of Proposition
\ref{urg}, namely to the Banach space $B_b,$ the operator $K$ defined
by $K f=\int_{{\mathcal S}_0} K(x,d y)f(y),$ and the cone $C$ of
nonnegative functions in $B_b.$ Note that $r_{_C}(K)$ coincides with
the spectral radius $r_{_K}$ in this case. It follows from
\textup{(\ref{comp1})} and the assumption $s(x,y) \in(c_1^{-1},c_1)$
that $r_{_K}>c^{-1/m}_1.$ Therefore it suffices to show that
\textup{(\ref{version})} implies $\rho_{_C}(K)=0.$ Since $\rho _{_C}(K)
\leq \rho_{_X}(K)$ (cf. \cite{nussbaum}, page~321), it is even
sufficient to show that \mbox{$\rho_{_X}(K)=0.$}\looseness=1

It will be convenient to use the notion of the Hausdorff measure of
noncompactness $\chi$ which is defined for a bounded subset $S$ of a
Banach space $X$ by
\[
\chi(S)=\inf\{d>0\dvtx  S \mbox{ has a finite } d\mbox{-net in } X\}.
\]
By finite $d$-net in $X$ we mean a finite subset $\{x_1,\ldots,x_n\}$
of $X$ such that for any $y \in{\mathcal S}$ there exists an index $j$
s.t. $\|y-x_j\|<d,$ where $\|\cdot\|$ is the norm
on~$X.$\vadjust{\goodbreak}

Let
\[
\chi(K):= \inf\{b \geq0\dvtx  \chi(K(S)) \leq b \chi(S) \mbox{ for
bounded subsets } S \mbox{ of } X\},
\]
and $\sigma(K):=\lim_{n \to\infty} \sqrt[n]{\chi(K^n)}.$ The Kuratowski
and Hausdorff measures of noncompactness are equivalent in the
following sense (cf. \cite{mnco}, page~4): $\chi(S) \leq\alpha(S) \leq2
\chi(S)$ for every bounded subset $S$ of $X.$ Thus, it suffices to show
that $\sigma(K)=0$ when \textup{(\ref{version})} holds. The latter
assertion follows from the following lemma.\looseness=1

\begin{lemma}
Let $X$ be a Banach space and $K$ be a bounded linear operator in $X.$
Further, let $\varepsilon>0$ be a positive constant and assume that
there is a compact operator $Q$ in $X$ such that $\|Q-K\| <
\varepsilon.$ Then, $\chi(K) \leq2 \varepsilon\|K\|.$
\end{lemma}

\begin{pf}
Fix a bounded set $S \subseteq X.$ Let $\{x_1,x_2,\ldots,x_n\}
\subseteq X$ be a finite $d$-net of $S$ for some $d>0.$ It suffices to
show that the set $K(S)$ has a finite $\eta_d$-net in $X,$ where we
denote $\eta_d:=2\varepsilon d\|K\|.$ Let $B_i,$ $i=1,\ldots,n,$ be the
balls in $X$ of radius $d$ and centered in $x_i.$ Then, $S
\subseteq\bigcup_{i=1}^n B_i$ and $K(S) \subseteq \bigcup_{i=1}^n
K(B_i).$ Therefore, it is sufficient to show that each set $K(B_i),$
$i=1,2,\ldots,n,$ has a finite $\eta_d$-net in~$X.$\looseness=1

Fix any $\delta>0.$ By the semi-homogeneity property of the measures of
noncompactness and their invariance under translations (cf.
\cite{mnco}, page~4) we can assume without loss of generality that
$d=1$ and consider only the unit ball $B_0$ centered at $0 \in X.$ Let
$Z:=\{z_1,z_2,\ldots,z_m\}$ be a finite $\delta$-net of the totally
bounded set $Q(B_0).$ Then, the balls of radius $\delta+\|K\|
\cdot\|K-Q\|$ with centers in $z_1,z_2,\ldots,z_m$ cover the set
$K(B_0).$ Indeed, for a point $x \in K(B_0),$ let $z(x) \in Z$ be such
that $\|Q x-z(x)\| \leq \delta.$ Then,
\begin{eqnarray*}
\|K x-z(x)\| &\leq& \|K x -Q x \| +\|Q x -z(x)\| \leq\|x\|
\cdot\|K-Q\|+\delta
\\
&\leq& \|K\| \cdot\|K-Q\|+\delta\leq\|K\| \cdot\|K-Q\|+\delta.
\end{eqnarray*}
This completes the proof of the lemma since $\delta>0$ is arbitrary.
\end{pf}
\end{appendix}

\section*{Acknowledgments}
I would like to thank my advisers at Technion-IIT (Haifa, Israel) Eddy
Mayer-Wolf and Ofer Zeitouni for numerous helpful conversations during
preparation of this work. I also thank Lev Grinberg for pointing out
references \cite{mnco,nussbaum} and the anonymous referee for his very
useful comments and suggestions.

\printaddresses


\begin{thebibliography}{99}
%Operator Theory: Advances and Applications, Vol. \textbf{55}.
%Translated from the 1986 Russian original by A. Iacob.
%b1 ###
\bibitem{mnco}
\textsc{Akhmerov, R. R.}, \textsc{Kamenski{\u\i}, M. I.}, \textsc
{Potapov, A. S.}, \textsc{Rodkina, A. E.} and \textsc{Sadovski{\u\i},
B. N.} (1992). \textit{Measures of Noncompactness and Condensing
Operators}.
Birkh\"auser, Basel.
\MR{1153247}\vadjust{\goodbreak}

%b2 ###
\bibitem{alsmeyer}
\textsc{Alsmeyer, G.} (1997). The {M}arkov renewal theorem and related
results. \textit{Markov Process. Related Fields} \textbf{3} 103--127.
\MR{1446921}

%Wiley Series in Probability and Mathematical Statistics: Applied
%Probability and Statistics. John
%b3 ###
\bibitem{asmus}
\textsc{Asmussen, S.} (1987). \textit{Applied Probability and Queues}.
Wiley, Chichester.
\MR{0889893}

%b4 ###
\bibitem{mac-dan}
\textsc{Athreya, K. B.}, \textsc{McDonald, D.} and \textsc{Ney, P.
{A}.} (1978). Limit theorems for semi-{M}arkov processes and renewal
theory for {M}arkov chains. \textit{Ann. Probab.} \textbf{6}
788--797.
\MR{0503952}

%b5 ###
\bibitem{athreya-ney}
\textsc{Athreya, K. B.} and \textsc{Ney, P. {A}.} (1978). A new
approach to the limit theory of recurrent {M}arkov chains.
\textit{Trans. Amer. Math. Soc.} \textbf{245} 493--501.
\MR{0511425}

%to the {P}erron--{F}robenius theory of nonnegative kernels on general
%state spaces. \textit{Math. Z.} \textbf{179} 507--529.

%b6 ###
\bibitem{brandt}
\textsc{Brandt, A.} (1986). The stochastic equation
$Y_{n+1}=A_nY_n+B_n$ with stationary coefficients. \textit{Adv. in
Appl. Probab.} \textbf{18} 211--220.
\MR{0827336}

%R. G.`Bartle.
%Pure and Applied Mathematics, Vol. \textbf{7}.
%Interscience Publishers, Inc., New York; Interscience Publishers,
%Ltd., London.

%b7 ###
\bibitem{saporta}
\textsc{de Saporta, B.} (2005). Tail of the stationary solution of the
stochastic equation {$Y_ {n+1}=a_ nY_ n+b_ n$} with {M}arkovian
coefficients. \textit{Stochastic Process. Appl.} \textbf{115}
1954--1978.
\MR{2178503}

%Examples}, 2nd ed. Duxbury Press, Belmont, CA.

%b8 ###
\bibitem{embre-goldie}
\textsc{Embrechts, P.} and \textsc{Goldie, C. M.} (1994). Perpetuities
and random equations. In \textit{Asymptotic Statistics}
(\textit{Prague}, \textit{1993}) 75--86. Contrib. Statist. Physica,
Heidelberg.
\MR{1311930}

%b9 ###
\bibitem{goldie}
\textsc{Goldie, C. M.} (1991). Implicit renewal theory and tails of
solutions of random equations. \textit{Ann. Appl. Probab.} \textbf{1}
126--166.
\MR{1097468}

%b10 ###
\bibitem{goldie-grubel}
\textsc{Goldie, C. M.} and \textsc{Gr{\"u}bel, R.} (1996). Perpetuities
with thin tails. \textit{Adv. in Appl. Probab.} \textbf{28} 463--480.
\MR{1387886}

%b11 ###
\bibitem{grey}
\textsc{Grey, D. R.} (1994). Regular variation in the tail behaviour of
solutions of random difference equations. \textit{Ann. Appl. Probab.}
\textbf{4} 169--183.
\MR{1258178}

%b12 ###
\bibitem{trakai75}
\textsc{Grincevi\v{c}jus, A. K.} (1975). On a limit distribution for a
random walk on the line. \textit{Lithuanian Math. J. } \textbf{15}
580--589.
\MR{0448571}

%b13 ###
\bibitem{trakai80}
\textsc{Grincevi\v{c}jus, A. K.} (1980). Products of random affine
transformations. \textit{Lithuanian Math. J.} \textbf{20} 279--282.
\MR{0605960}

%b14 ###
\bibitem{kesten-randeq}
\textsc{Kesten, H.} (1973). Random difference equations and renewal
theory for products of random matrices. \textit{Acta Math.}
\textbf{131} 207--248.
\MR{0440724}

%b15 ###
\bibitem{kesten-rem}
\textsc{Kesten, H.} (1974). Renewal theory for functionals of a
{M}arkov chain with general state space. \textit{Ann. Probab.}
\textbf{2} 355--386.
\MR{0365740}

%additifs en \'economie stochastique. \textit{C. R. Acad. Sci. Paris
%S\'er. A} \textbf{279} 33--36.

%b16 ###
\bibitem{mars}
\textsc{Mayer-Wolf, E.}, \textsc{Roitershtein, A.} and \textsc
{Zeitouni, O.} (2004). Limit theorems for one-dimensional transient
random walks in {M}arkov environments. \textit{Ann. Inst. H. Poincar\'e
Probab. Statist.} \textbf{40} 635--659.
\MR{2086017}

%Communications and Control Engineering Series.
%b17 ###
\bibitem{meyn-tweedie}
\textsc{Meyn, S. P.} and \textsc{Tweedie, R. L.} (1993). \textit{Markov
Chains and Stochastic Stability}.
Springer, London.
\MR{1287609}

%Cambridge Tracts in Mathematics, Vol. \textbf{83}.
%b18 ###
\bibitem{nummelin}
\textsc{Nummelin, E.} (1984). \textit{General Irreducible {M}arkov
Chains and Nonnegative Operators}.
Cambridge Univ. Press.
\MR{0776608}

%b19 ###
\bibitem{nussbaum}
\textsc{Nussbaum, R. D.} (1981). Eigenvectors of nonlinear positive
operators and the linear {K}rein--{R}utman theorem. \textit{Fixed Point
Theory}. %(\textit{Sherbrooke}, \textit{Que.}, \textit{1980})
\textit{Lecture Notes in Math.} \textbf{886} 309--330. Springer, Berlin.
\MR{0643014}

%b20 ###
\bibitem{samorachev}
\textsc{Rachev, S. T.} and \textsc{Samorodnitsky, G.} (1995). Limit
laws for a stochastic process and random recursion arising in
probabilistic modelling. \textit{Adv. in Appl. Probab.} \textbf{27}
185--202.
\MR{1315585}

%b21 ###
\bibitem{revuz}
\textsc{Revuz, D.} (1975). \textit{Markov Chains}.
North-Holland, Amsterdam.
\MR{0415773}

%b22 ###
\bibitem{shur}
\textsc{Shurenkov, V. M.} (1984). On {M}arkov renewal theory.
\textit{Teor. Veroyatnost. i Primenen.} \textbf{29} 248--263.
\MR{0749913}

%b23 ###
\bibitem{vervaat}
\textsc{Vervaat, W.} (1979). On a stochastic difference equation and a
representation of nonnegative infinitely divisible random variables.
\textit{Adv. in Appl. Probab.} \textbf{11} 750--783.
\MR{0544194}\

%b24 ###
\bibitem{yosida}
\textsc{Yosida, K.} and \textsc{Kakutani, S.} (1941).
Operator-theoretical treatment of {M}arkoff's process and mean ergodic
theorem. \textit{Ann. of Math.} (\textit{2}) \textbf{42}
188--228.
\MR{0003512}

\end{thebibliography}
\end{document}